\documentclass[review,3p]{elsarticle}

\usepackage{nameref}

\usepackage{graphicx}
\usepackage{enumitem}
\usepackage{amsmath}
\usepackage{bm}
\usepackage{makecell}
\usepackage{caption, subcaption}
\usepackage{xcolor}
\usepackage{algpseudocode}
\usepackage{algorithm}
\usepackage{standalone}
\usepackage{tikz,pgf}
\usepackage{amsfonts}
\usepackage{multirow}
\usepackage{textcomp}
\usepackage{placeins,float}

\usepackage{titlesec}
\titleformat{\paragraph}[runin]{\normalfont\normalsize\bfseries}{\theparagraph}{1em}{}
\titlespacing*{\paragraph}{0pt}{0pt}{1em}{}

\usepackage{lipsum}
\makeatletter \def\ps@pprintTitle{\let\@oddhead\@empty \let\@evenhead\@empty \def\@oddfoot{\hfill\today} \let\@evenfoot\@oddfoot} \makeatother

\begin{document}

\begin{frontmatter}

\title{Branch-and-bound algorithm for efficient reliability analysis \\ of general coherent systems}

\author[1]{Ji-Eun Byun\corref{cor1}}
\ead{Ji-Eun.Byun@glasgow.ac.uk}
\address[1]{James Watt School of Engineering, University of Glasgow, Glasgow, United Kingdom}

\author[2]{Hyeuk Ryu}
\ead{hyeuk.ryu@ga.gov.au}
\address[2]{Geoscience Australia, Australia}

\author[3]{Daniel Straub}
\ead{straub@tum.de}
\address[3]{Engineering Risk Analysis Group, Technical University of Munich, Germany}

\cortext[cor1]{Corresponding author}

\begin{abstract}
Branch and bound algorithms have been developed for reliability analysis of coherent systems. They exhibit a set of advantages; in particular, they can find a computationally efficient representation of a system failure or survival event, which can be re-used when the input probability distributions change over time or when new data is available. However, existing branch-and-bound algorithms can handle only a limited set of system performance functions, mostly network connectivity and maximum flow. Furthermore, they run redundant analyses on component vector states whose system state can be inferred from previous analysis results. This study addresses these limitations by proposing \textit{\textbf{b}ranch and bound for \textbf{r}eliability analysis of general \textbf{c}oherent systems} (BRC) algorithm: an algorithm that automatically finds minimal representations of failure/survival events of \textit{general} coherent systems. Computational efficiency is attained by dynamically inferring importance of component events from hitherto obtained results. We demonstrate advantages of the BRC method as a real-time risk management tool by application to the Eastern Massachusetts highway benchmark network.
\end{abstract}

\begin{keyword}
    system reliability \sep large-scale systems \sep discrete-state components \sep non-simulation-based approach \sep branch-and-bound algorithm \sep reliability mapping \sep real-time model updating
\end{keyword}

\end{frontmatter}

\section{Introduction}

Reliability analysis of large-scale systems remains challenging because of the high-dimensional probability distribution of a large number of components and the high non-linearity of failure domains. Those systems often modelled with discrete-state components to represent discrete degradation of functionality (e.g. traffic capacity of roadways conditioned on their discrete damage states \cite{LeeSonGarLim11,BoaGuiGarMur22,NiuXiaXu24}). As a result, computational complexity increases exponentially with the number of components since there are an exponential number of possible combinations of component states. In addition, with discrete-state components, there are often a greater number of disconnected failure domains than with continuous-state ones. Missing important failure domains is bound to underestimate system failure probability, which is undesirable from the perspective of risk management. Discrete-state components also lead to non-smooth surface of a system performance function (a.k.a. a limit-state function), making it challenging to efficiently explore the surface. 

Simulation methods have been developed to handle a large number of components such as importance sampling \cite{DasJoh24, KanPapStr23} and subset simulation \cite{CheWanBro22, DuLuoWan19}. However, those methods were proposed to address continuous-state components, and thus often become inefficient with discrete-state components. To address this issue, various advanced simulation methods have been enhanced towards applications to systems with discrete inputs, such as variants of subset simulation \cite{LeeWangSong23, ChanPapaStra22,YuHuanWen21} or cross-entropy-based importance sampling \cite{ChaPapStr24}. 
\cite{LeeByuSad23, ZwiChaPap24} proposed combining a sampling method with a Bayesian network (BN)-based representation.
Despite improved performance, these methods still need a significant number of system function evaluations, in particular for systems with high reliability and large number of components. 

Another popular approach is using surrogate models of a system performance function such as Gaussian process \cite{KimYiSon24, MouSud23} and polynomial chaos expansion \cite{LiuEtal21, SchPapEhr20,WenLiuZhaZha24}. However, these methods also assume smooth surface of a system function, which is not a valid assumption given discrete-state components. Deep neural networks (DNNs) present a promising potential for general input and output types. DNNs have been applied for system reliability analyses such as dynamic response of structural systems \cite{KimKwoSon24}, regional loss assessment by damaged buildings \cite{KanKimKwo23}, network connectivity \cite{DavYodLeYad23,ShiBehZhoHu24,HuaLin24}, and network maximum flow \cite{HuaHuaLin22,HuaHuaLin23}. However, surrogate-based approaches still need a relatively large number of training data sets (e.g. \cite{HuaHuaLin22} requires 10,000 training data for a system with 20 components). This is because these methods make minimal assumptions on a given system event. When useful information on system characteristics is available, they may not be the best solution. 

Non-simulation-based methods have been continuously studied to solve network and general system reliability problems. Among the most popular approaches are branch-and-bound algorithms, which efficiently decompose a system event space into survival and failure domains \cite{LiHe02, JanLai08}. Their strategies are largely two-fold. The first strategy identifies all minimal conditions of failure and/or survival (often referred to as minimal cuts and minimal paths) and then proceeds to decomposition of a system event space \cite{Yeh21,LinChe17, ForKagMah19, NiuWanXu20,ZhoBaiTaoXu23,Cha24}. The second strategy is to identify survival and failure conditions and decompose an event space, at the same time \cite{LiHe02, LimSon12,Koz24,JanLai08}. The first strategy is more efficient when all minimal conditions can be identified. However, as the number of components increases, the number of those conditions often becomes too large for exhaustive identification. The second strategy has an advantage that an algorithm can terminate without identifying all conditions, in which case it can provide a bound on system failure probability rather than an exact value. Still, a too large number of components often makes the bound too wide to be useful.

Branch-and-bound algorithms exploit \textit{coherency} in the system function. When a system is coherent, its state does not worsen by improving component states. This is crucial information: once some component states are known to lead to system failure, one can instantly infer that all worse states lead to system failure without running further system simulation. Consequently, these methods often require a much smaller number of system function evaluations (often much less than 1,000). Despite such efficiency, previous branch-and-bound algorithms are limited in that they were developed focusing on specific system performance functions, mostly network connectivity and maximum flow. As a result, their applicability has been significantly hampered as a new algorithm needs to be developed to handle even a slight modification of a system performance definition.

The objectives of the present study are two-fold: (1) handling general coherent systems and (2) removing limits in the addressable number of components. To this end, this study proposes a new algorithm namely \textit{\textbf{b}ranch and bound for \textbf{r}eliability analysis of general \textbf{c}oherent systems} (BRC) algorithm. The algorithm consists of a set of effective strategies to identify minimal representations of failure and survival component states that decide system failure and survival. To mitigate the algorithm's limitations on the number of components, we derive formulas for hybrid inference that combine branch-and-bound with Monte Carlo Simulation (MCS). We compare the algorithm with a previous branch-and-bound algorithm of \cite{JanLai08} to demonstrate its improved efficiency despite its general applicability. We also apply the BRC algorithm to Eastern Massachusetts (EMA) benchmark network, which demonstrates its applicability to general and large-scale systems. 

The article is organised as follows. In Section~\ref{sec:setting}, we introduce notations and the scope of the work. Section~\ref{sec:brc} introduces the proposed algorithm and illustrates it with a basic example. Section~\ref{sec:ex} presents two numerical examples. In Section~\ref{sec:dis}, we provide detailed discussions on the BRC algorithm. Finally, Section~\ref{sec:con} ends the paper with concluding remarks.

\section{Notations and scope}\label{sec:setting}

In this paper, a random variable and an assignment of its state are respectively denoted by upper and lower cases, e.g. $X$ and $x$. When referring to a set of random variables, we use bold symbols, e.g. $\bm{X}$ and $\bm{x}$. For simplicity, an assignment $X=k$ is abbreviated as $x^k$. Adopting the terminology of Koller and Friedman \cite{KolFri09}, the set of states that a random variable $X$ can take is denoted $Val(X)$. Given a state $\bm{x}$, we refer to the set of random variables that a state $\bm{x}$ is defined over as $Scope[\bm{x}]$. The state of $X\in Scope[\bm{x}]$ in $\bm{x}$ is denoted as $\bm{x}\langle X \rangle$.

We consider a \textit{system} consisting of $N$ \textit{components} that have discrete states of functionality. The states of a system and its components are represented by random variables $S$ and $\bm{X} = \{X_1, \cdots, X_N\}$. Each component can take one of $K$ discrete states $0,1,\ldots,K-1$, i.e. $Val(X_n)=\{0,\ldots,K-1\}$. Without loss of generality, it is assumed that a higher state $x_n$ represents a corresponding component's better functionality (e.g. greater flow capacity and shorter travel time). On the other hand, $S$ takes a binary state, 0 for system failure and 1 for system survival, i.e. $Val(S)=\{0,1\}$. The state of $S$ is evaluated by a \textit{system performance function}, denoted as $\Phi(\cdot)$, which is a deterministic function that accepts as input a \textit{component vector state} (i.e. an assignment on $\bm{X}$) and returns a corresponding state of $S$, i.e. $s=\Phi(\bm{x})$.

\begin{figure}[h!]
\centering
    \begin{subfigure}[b]{0.49\textwidth}
        \centering
        \newcommand\xdist{1.5}

\begin{tikzpicture}
    \node[circle, draw, fill=gray!30] (1) at (0,0) {$n_1$};
    \node[circle, draw] (2) at (\xdist,0) {$n_2$};
    \node[circle, draw, fill=gray!30] (3) at (2*\xdist,0) {$n_3$};
    
    \draw[-] (1) edge node[above] {$e_1$} (2);
    \draw[-, bend left] (2) edge node[above] {$e_2$} (3);
    \draw[-, bend right] (2) edge node[below] {$e_3$} (3);
\end{tikzpicture}
        \caption{Graph and three fragile edges.}
        \label{subfig:toy}
    \end{subfigure}
    \hfill
    \begin{subfigure}[b]{0.49\textwidth}
        \centering
        \newcommand\xdist{1.5}

\begin{tikzpicture}
    \node[circle, draw, fill=gray!30] (1) at (0,0) {$n_1$};
    \node[circle, draw] (2) at (\xdist,0) {$n_2$};
    \node[circle, draw, fill=gray!30] (3) at (2*\xdist,0) {$n_3$};
    
    \draw[-] (1) edge node[above] {$e_1$} (2);
    \draw[dashed, bend left] (2) edge node[above] {$e_2$} (3);
    \draw[dashed, bend right] (2) edge node[below] {$e_3$} (3);
\end{tikzpicture}
        \caption{A failure mode where $e_2$ and $e_3$ fail together.}
        \label{subfig:toy2}
    \end{subfigure}
    
    \begin{subfigure}[b]{0.49\textwidth}
        \centering
        \tikzstyle{block} = [draw, rectangle, minimum height=0.8cm, minimum width=1.4cm]

\begin{tikzpicture}[thick, main/.style = {draw, circle}, node distance=1.5cm] 

\newcommand\xdist{1}
\newcommand\ydist{1.7}
\newcommand\xtext{2.7}
\newcommand\ytext{0.3}
\newcommand\ytexti{1.2}
\newcommand\ys{1cm} 
\newcommand\nsz{0.75cm}

\node[main, minimum size = \nsz] at (0,0) (b0) {$b_0$};

\node[main] at (-\xdist,-\ydist) (b1) {$b_1$};
\draw[->] (b0) -- (b1) node[midway,above left] {$X_2 \leq 0$};

\node[main] at (\xdist,-\ydist) (b2) {$b_2$};
\draw[->] (b0) -- (b2) node[midway,above right] {$X_2 \geq 1$};

\node[main, fill=black!10] at (-\xdist-\xdist,-\ydist-\ydist) (b3) {$b_3$};
\draw[->] (b1) -- (b3) node[midway,above left] {$X_3 \leq 0$};

\node[main] at (0,-\ydist-\ydist) (b4) {$b_4$};
\draw[->] (b1) -- (b4) node[midway,above right] {$X_3 \geq 1$};

\end{tikzpicture}
        \caption{An example decomposition procedure. Given the failure mode in Fig.~\ref{subfig:toy2}, one can infer that all component vector states in branch $b_3$ cause system failure.}
        \label{subfig:toy3}
    \end{subfigure}

    \caption{Illustrative example}\label{fig:toy}
\end{figure}
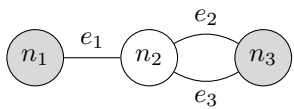
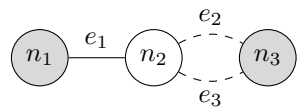
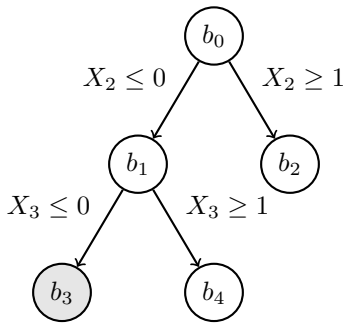

For example, consider the network in Figure~\ref{subfig:toy}. There are three edges subject to failure, $e_1$, $e_2$, and $e_3$, which are represented by random variables $X_1$, $X_2$, and $X_3$. Each $X_n$, $n=1,2,3$, is binary (i.e. $Val(X_n)=\{0,1\}$), where 0 for failure and 1 for survival. System failure is defined as nodes $n_1$ and $n_3$ being disconnected, i.e. the system performance function returns 0 if nodes 1 and 3 are not connected and 1 otherwise.  

The scope of this paper is confined to \textit{coherent} systems, i.e. an increase (decrease) in a component vector state does not deteriorate a system state from survival to failure (improve from failure to survival). Coherency is commonly observed in real-world situations such as a transportation network modelled by maximum flow or shortest path\footnote{Incoherency can arise from suboptimal decisions by operators and/or users. For instance, a power grid is an incoherent system if it is subject to cascading failures by excessive power flows from a distribution node. Such overflows occur when a substation is unaware of the failure of downstream components and does not adjust its supply. Another example is a transportation system in which users make decisions that are suboptimal with respect to overall system performance.}.

In this paper, coherency is exploited for computational efficiency by the fact that if a state  $\bm{r}$ with $Scope[\bm{r}]\subseteq \bm{X}$ is known to cause a system failure (survival), it can represent an upper (lower) bound for the system state. That is, the same system state is guaranteed if a state $\bm{x}$ satisfies the upper (lower) bounds set by $\bm{r}$, i.e. $\bm{x}\langle X \rangle \leq \bm{r}\langle X \rangle$ ($\bm{x}\langle X \rangle \geq \bm{r}\langle X \rangle$) for $\forall X \in Scope[\bm{r}]$. In addition, as long as the bounds $\bm{r}$ are satisfied, the components $X \notin Scope[\bm{r}]$ do not affect the system state. For instance, Fig.~\ref{subfig:toy2} illustrates a failure mode where $e_2$ and $e_3$ fail together, i.e. $\bm{r}=(x_2^0,x_3^0)$ leads to $s^0$. From this result, one can infer, without running the system performance function, that $(x_1^0,x_2^0,x_3^0)$ and $(x_1^1,x_2^0,x_3^0)$, which satisf $X_2 \leq 0$ and $X_3 \leq 0$, lead to system failure. 

To formalise this concept, we introduce {\it rules}. A rule $\gamma=(\bm{r}, s)$ is defined by a condition state $\bm{r}$ and an associated system state $s$. If $s=0$, a rule is called a \textit{failure} rule, and if $s=1$, a \textit{survival} rule. We denote the set of rules as $\mathcal{R}$, that of failure rules as $\mathcal{R}_f$, and that of survival rules as $\mathcal{R}_s$ (i.e. $\mathcal{R}=\mathcal{R}_f \cup \mathcal{R}_s$). The failure mode illustrated in Fig.~\ref{subfig:toy2} can be represented by a failure rule $\gamma=\left((x_2^0,x_3^0), 0\right)$. A rule is called \textit{minimal} if removal of any element of $\bm{r}$ no longer guarantees an associated $s$. For instance, the previous failure rule $\gamma=\left((x_2^0,x_3^0), 0\right)$ is minimal since removing either $x_2^0$ or $x_3^0$ no longer guarantees system failure. A failure (survival) rule can be considered a generalised form of \textit{cut-sets} (\textit{link-sets}) for multi-state components. While cut-sets and link-sets are fully defined by their associated binary components, rules require associated states as well to handle multi-state components.

Given a complete set of rules (i.e. one with which all existing component vector states can be classified to either system failure and survival), the system failure probability is equivalent to the joint probability of the rules, i.e.
\begin{equation}\label{eq:pf_rule}
    \begin{aligned}
        P( s^0 ) & = P\left( \cup_{(\bm{r},0)\in \mathcal{R}_f} \left\{ \cap_{X \in Scope[\bm{r}]} X \leq \bm{r}\langle X \rangle \right\} \right) \\
        & = 1 - P\left( \cup_{(\bm{r},1)\in \mathcal{R}_s} \left\{ \cap_{X \in Scope[\bm{r}]} X \geq \bm{r}\langle X \rangle \right\} \right).
    \end{aligned}
\end{equation}
However, the computational cost for evaluating such joint probabilities increases exponentially with the number of rules, hence evaluating \eqref{eq:pf_rule} becomes infeasible when dealing with a large number of components. Therefore, one often resorts to branch-and-bound approaches that decompose the event space into disjoint \textit{branches} \cite{LiHe02,JanLai08}.

We denote a branch as $b=\left( \bm{l}, \bm{u}, s_l, s_u, p \right)$, where $\bm{l}$ and $\bm{u}$ stand for lower and upper bounds on the component vector state, and $s_l$ and $s_u$ represent the system state corresponding to $\bm{l}$ and $\bm{u}$, respectively. $s_l$ and $s_u$ can take three values, 0 for system failure, 1 for system survival, and u for unknown system state. $p$ indicates the probability of a branch, defined as
\begin{equation}\label{eq:braPro}
    p=P(b) = P(\bm{l} \leq \bm{X} \leq \bm{u}).       
\end{equation}
If calculating the probability $p$ is computationally expensive, which is often the case when component events are dependent, one can approximate it by
\begin{equation}\label{eq:braProApr}
p \approx \Tilde{P}\left( b \right) =
        \prod_{X\in \bm{X}} P\left( \bm{l} \langle X \rangle \leq X \leq \bm{u} \langle X \rangle \right),
\end{equation}
which assumes statistical independence between component variables. This assumption should work effectively for the purpose of the proposed algorithm, with the exception of some extreme cases, e.g. very high correlation between component events.

For example, in Fig.~\ref{subfig:toy3}, suppose that the components are independent with failure probabilities $P(x_1^0)=0.1$, $P(x_2^0)=0.2$, and $P(x_3^0)=0.3$. Then, the branch $b_3$ is represented as $b_3=\left( (0,0,0), (1,0,0), 0, 0, 0.06 \right)$, where the probability is calculated as $P(X_2 \leq 0) \cdot P(X_3 \leq 0) = 0.2 \cdot 0.3 = 0.06$. Using the failure rule $\gamma=\left((x_2^0,x_3^0), 0\right)$, one can infer that both the upper and the lower bounds of $b_3$ lead to system failure. By contrast, branch $b_2$ is $b_2=\left( (0,1,0), (1,1,1), \text{u}, \text{u}, 0.8 \right)$. Their associated system states remain unknown since the only known failure rule $\gamma$ is not satisfied by both the upper and lower bound of $b_2$. The probability of $b_2$ is calculated as $P(X_2\geq1) = 0.8$.

A branch is considered {\it specified} if $s_l = s_u \neq \text{u}$, i.e. all included component vector states are known to lead to the same system state; if not, a branch is considered {\it unspecified}. We refer to branches with $s_l=s_u=0$ as a \textit{failure} branches and those with $s_l=s_u=1$ as a \textit{survival} branches. We denote the set of branches as $\mathcal{B}$, that of failure branches as $\mathcal{B}_f$, that of survival branches as $\mathcal{B}_s$, and that of unspecified branches $\mathcal{B}_u$, where $\mathcal{B}=\mathcal{B}_f \cup \mathcal{B}_s \cup \mathcal{B}_u$. For instance, in Fig.~\ref{subfig:toy3}, $b_3$ is a failure branch, while $b_2$ and $b_4$ are unspecified branches.

Once branches are obtained using rules, since they are disjoint to each other, a system failure probability can be calculated simply as a sum of the probabilities of failure branches. Given that there are no unspecified branches, \eqref{eq:pf_rule} becomes
\begin{equation}
    P( s^0 ) = \sum_{b\in \mathcal{B}_f} P(b) = 1 - \sum_{b\in \mathcal{B}_s} P(b),
\end{equation}
where the computational cost is linear in the number of branches. The number of branches can increase rapidly, sometimes exponentially, with the number of rules, but this increase can be effectively mitigated by optimising the decomposition sequence (for which we propose strategies in Section~\ref{sec:dec}). 

\section{Branch-and-bound for reliability analysis of coherent systems (BRC) algorithm}\label{sec:brc}

\subsection{The overall structure}\label{sec:brc_overview}

\begin{figure}
    \centering
    \includegraphics[width=0.9\textwidth]{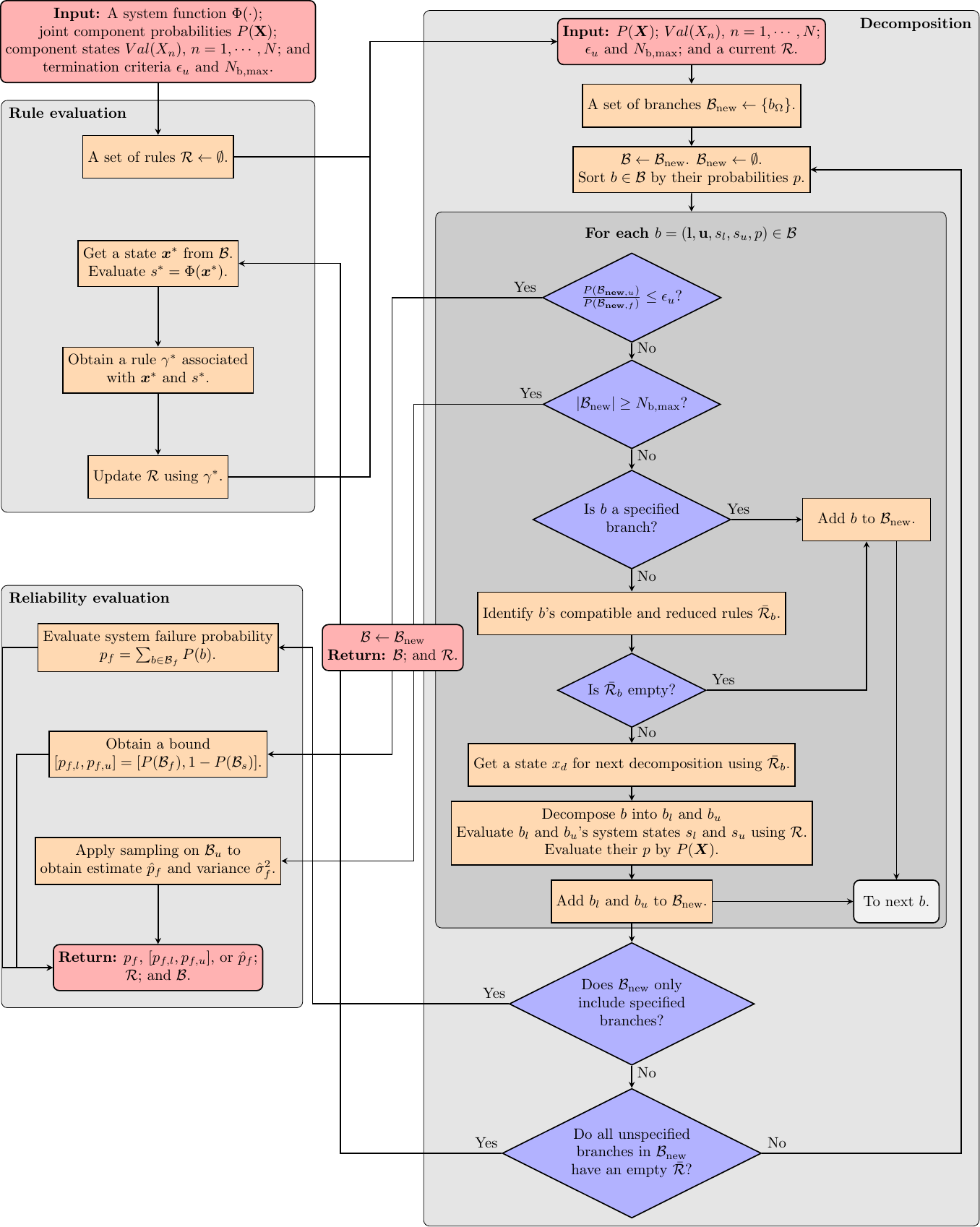}
    \caption{The BRC algorithm consisting of rule evaluation (upper left-hand side), decomposition (right-hand side), and reliability evaluation (lower left-hand side).}
    \label{fig:brc}
\end{figure}

The BRC algorithm is summarised in Fig.~\ref{fig:brc}. The algorithm consists of three modules: rule evaluation, decomposition, and reliability evaluation. The rule evaluation module runs a system performance function and thereby finds new rules. The decomposition module is run every instance a new rule is found in the rule evaluation procedure. It identifies an optimal decomposition given a current set of rules and informs on a next step of the rule evaluation. Finally, the reliability evaluation module concludes the algorithm by evaluating a system failure probability. 
We present the decomposition module in Section~\ref{sec:dec}, the rule evaluation process in Section~\ref{sec:rule_eval}, and the final reliability evaluation in Section~\ref{sec:stop}. In Section~\ref{sec:toy}, we illustrate the application of the BRC algorithm using the example network of Fig.~\ref{fig:toy}. 

\subsection{Decomposition}\label{sec:dec}

\subsubsection{Inputs/outputs, strategies, and initialisation}

The decomposition algorithm performs an iterative depth-first decomposition until no further specified branches can be produced using a current $\mathcal{R}$. The algorithm takes as inputs a current set of rules $\mathcal{R}$ together with $P(\bm{X})$ , $Val(X_n)$, $n=1,\cdots,N$, and termination parameters $\epsilon_u$ and $N_{\text{b,max}}$. The algorithm has three purposes. First, it enables early termination when a given system scale appears to be so large that it needs hybrid inference. Second, it heuristically optimises the decomposition sequence to minimise the number of branches resulting from the given rules. Third, the decomposition result, i.e. a set of branches $\mathcal{B}$, is exploited during the rule evaluation process for deciding for which component vector state a system performance function evaluation should be performed next. 

We introduce three strategies to make the algorithm efficient for these tasks. First, the algorithm considers branches with higher probabilities first. This intends to expedite narrowing down a bound on system failure probability. Second, it obtains \textit{reduced, compatible} rules to evaluate the relative importance of component states from a branch-specific, not global, perspective. Third, to minimise the number of system function runs, the algorithm does not use a system performance function $\Phi(\cdot)$ to decide the system states of branches. Rather, it compares the lower and upper bounds of a branch with given rules in $\mathcal{R}$. These strategies are explained in detail in the following subsections.

For initialisation, the current set of branches is set to $\mathcal{B}_{\text{new}}=\{b_{\Omega}\}$. $b_{\Omega}$ represents the total event space, $b_{\Omega}=(\bm{x}_{\mathrm{min}}, \bm{x}_{\mathrm{max}}, s_{l,\mathrm{min}}, s_{u,\mathrm{max}}, 1)$, where $\bm{x}_{\mathrm{min}}$ and $\bm{x}_{\mathrm{max}}$ are respectively the worst and best states of $\bm{X}$. $s_{l,\mathrm{min}}$ and $s_{u,\mathrm{max}}$ are corresponding states to $\bm{x}_{\mathrm{min}}$ and $\bm{x}_{\mathrm{max}}$, respectively. The branch by definition has probability $1$.

After initialisation, a for-loop (the dark box in Fig.~\ref{fig:brc}) is iterated to decompose branches in $\mathcal{B}$ and thereby identify more specified branches. Before going into the for-loop, branches are sorted by their probabilities so that branches with higher probabilities can be considered first.

\subsubsection{Early termination}

Once the algorithm enters the for-loop, there are two check points for early termination of the BRC algorithm (i.e. termination with unspecified branches remaining in $\mathcal{B}_{\text{new}}$). First, the algorithm is terminated if the ratio between the probabilities of unknown and failure, i.e. $P(\mathcal{B}_{\text{new},u}) / P(\mathcal{B}_{\text{new},f})$, is less than a termination threshold $\epsilon_u$ (e.g. 0.05). The value represents the ratio of the final bound width to the lower bound of a system failure probability estimate. When one sets $\epsilon_u = 0$, the algorithm returns an exact value. The second exit point is when the number of branches $\lvert \mathcal{B}_{\text{new}} \rvert$ exceeds a user-specified limit $N_{\text{b,max}}$ (e.g. 50,000). In this case, the BRC algorithm stops decomposition and moves on to hybrid inference by sampling on the unspecified branches $\mathcal{B}_u$ (c.f. Section~\ref{sec:stop}).

\subsubsection{Deciding on the decomposition of a branch}\label{sec:decomp}

Given a current branch $b$ in the for-loop, if neither of the early termination criteria is satisfied, the algorithm checks if $b$ is available for decomposition. There are two points to check. First, if $b$ is a specified branch, i.e. $s_l=s_u\neq \text{u}$, it does not need decomposition. Otherwise, the algorithm evaluates $b$'s set of reduced, compatible rules $\Bar{\mathcal{R}}$, which are defined in Section~\ref{sec:bar_r}. When $\Bar{\mathcal{R}}$ is empty, no further decomposition is possible using a current $\mathcal{R}$.

When it is decided that $b$ will not be further decomposed, it is added to $\mathcal{B}_{\text{new}}$, and the for-loop moves on to the next $b$. Otherwise, the algorithm proceeds to the decomposition of $b$ (c.f. Section~\ref{sec:xd}).

\subsubsection{Obtaining a set of compatible and reduced rules}\label{sec:bar_r}

Given an unspecified branch $b$, the algorithm evaluates a set of \textit{compatible} and \textit{reduced} rules, $\Bar{\mathcal{R}}_b$, which are used to check whether further decomposition is possible (c.f. Section~\ref{sec:decomp}) and to decide the next component state for decomposition (c.f. Section~\ref{sec:xd}). $\Bar{\mathcal{R}}_b$ is defined as follows.

Given a branch $b=\left( \bm{l}, \bm{u}, s_l, s_u, p \right)$, we define \textit{compatible} rules as those whose condition can be satisfied by a subset of component vector states within $b$. In other words, a rule $\gamma=(\bm{r}, s)$ is compatible with $b$ if
\begin{subequations}
\begin{equation}\label{eq:com_rf}
    \bm{r} \langle X \rangle \geq \bm{l} \langle X \rangle, \forall X \in Scope[\gamma] \text{\quad for $s=0$ and}
\end{equation}
\begin{equation}\label{eq:com_rs}
    \bm{r} \langle X \rangle \leq \bm{u} \langle X \rangle, \forall X \in Scope[\gamma] \text{\quad for $s=1$.}
\end{equation}
\end{subequations}
In \eqref{eq:com_rf} where $\gamma$ is a failure rule, the condition $\bm{r}$ ensures that the branch can still satisfy the rule. \eqref{eq:com_rs} applies the same logic to a survival rule.

A rule $\gamma$ compatible to $b$ can be \textit{reduced} to include only the conditions that are not already satisfied by $b$. That is, a reduced rule $\Bar{\gamma}=(\Bar{\bm{r}},s)$ inherits the condition states of $\gamma$ but has a reduced scope as
\begin{subequations}
    \begin{equation}
        \Bar{\bm{r}} = \{\bm{r}\langle X \rangle: X\in Scope[\gamma] \text{ and } \bm{r}\langle X \rangle < \bm{u} \langle X \rangle\} \text{\quad if $s=0$  and}
    \end{equation}
    \begin{equation}
        \Bar{\bm{r}} = \{\bm{r}\langle X \rangle: X\in Scope[\gamma] \text{ and } \bm{r}\langle X \rangle > \bm{l} \langle X \rangle \} \text{\quad if $s=1$.}
    \end{equation}
\end{subequations}

For instance, in the system of Fig.~\ref{subfig:toy}, suppose the current rule set is 
$\mathcal{R}=\{\gamma_1,\gamma_2,\gamma_3\}=\{\left((x_1^1,x_2^1), 1\right), \\\left((x_1^1,x_3^1), 1\right),\left((x_2^0,x_3^0), 0\right)\}$, and the branch is $b_1=(\bm{l}_1, \bm{u}_1, s_{l,1}, s_{u,1},p_1)=\left((0,0,0),(1,1,0),0,1,0.3\right)$. In this case, the compatible and reduced rules of $b_1$ are $\Bar{\mathcal{R}}_1=\{\gamma_1, \Bar{\gamma}_3\}$: $\gamma_2$ is excluded since $\bm{r}_2\langle X_3 \rangle = 1$ is incompatible with $\bm{u}_1\langle X_3 \rangle = 0$. $\gamma_3$ is reduced to $\Bar{\gamma}_3=\left((x_2^0), 0\right)$ since the condition $\bm{r}_3\langle X_3\rangle =0$ is always satisfied by $\bm{u}_1\langle X_3 \rangle = 0$. 

\subsubsection{Selecting a component state for decomposition}\label{sec:xd}

A non-empty $\Bar{\mathcal{R}}_b$ is used to select the component state $x_d$ for decomposing $b$. To minimise the number of branches, we propose selecting $x_d$ as follows. First, one orders component variables $\bm{X}$ as $\bm{X}'$ by their number of appearances in $\Bar{\mathcal{R}}_b$. Then, rules $\Bar{\mathcal{R}}_b$ are ordered by their probabilities conditioned on $b$ and stored in $\Bar{\mathcal{R}}_b'$. The conditional probability of a rule $\Bar{\gamma}=(\Bar{\bm{r}},s)\in \Bar{\mathcal{R}}_b$ is
\begin{equation}\label{eq:conPro}
P\left( \Bar{\gamma} | b \right) \propto
    \begin{cases}
        P\left( \cap_{X\in Scope[\Bar{\gamma}]} \left\{ \bm{d} \langle X \rangle \leq X \leq \Bar{\bm{r}} \langle X \rangle \right\} \right) & \text{for $s=0$,} \\
        P\left( \cap_{X\in Scope[\Bar{\gamma}]} \left\{ \Bar{\bm{r}} \langle X \rangle \leq X \leq \bm{u} \langle X \rangle \right\} \right) & \text{for $s=1$.}        
    \end{cases}
\end{equation}
The proportionality constant $1/P(b)$ in Eq. \eqref{eq:conPro} is ignored as it is the same for all rules $\Bar{\gamma}=\Bar{\mathcal{R}}_b$. If \eqref{eq:conPro} is too expensive to calculate, the joint probability can be approximated similarly to \eqref{eq:braProApr} using marginal distributions, i.e.
\begin{equation}
    \Tilde{P}\left( \Bar{\gamma} | b \right) \propto
    \begin{cases}
        \prod_{X\in Scope[\Bar{\gamma}]} P\left( \bm{d} \langle X \rangle \leq X \leq \Bar{\bm{r}} \langle X \rangle \right) & \text{for $s=0$,} \\
        \prod_{X\in Scope[\Bar{\gamma}]} P\left( \Bar{\bm{r}} \langle X \rangle \leq X \leq \bm{u} \langle X \rangle \right) & \text{for $s=1$,}
        
    \end{cases}
\end{equation}
The reason because we propose performing all these evaluations using $\Bar{\mathcal{R}}_b$, not the full set of rules $\mathcal{R}$, is to ensure that importance of components and rules can be evaluated conditional on the specific positions of a given branch. 

Running through the ordered sets $\bm{X}'$ and $\Bar{\mathcal{R}}_b'$, we propose selecting the first component state that allows for decomposition. In other words, one selects the first state $\Bar{\bm{r}} \langle X \rangle$ that satisfies
\begin{subequations}
\begin{equation}\label{eq:com_rf2}
    \bm{l}\langle X \rangle \leq \Bar{\bm{r}} \langle X \rangle < \bm{u}\langle X \rangle \text{\quad for $s=0$ and}
\end{equation}
\begin{equation}\label{eq:com_rs2}
    \bm{l}\langle X \rangle < \Bar{\bm{r}} \langle X \rangle \leq \bm{u}\langle X \rangle \text{\quad for $s=1$.}
\end{equation}
\end{subequations}
If these conditions are not satisfied, decomposition cannot be performed as a given state is out of the bound of a given $b$. Once $\Bar{\bm{r}} \langle X \rangle$ is selected, $x_d$ is set to
\begin{equation}
    x_d = 
    \begin{cases}
        x_d = \Bar{\bm{r}} \langle X \rangle + 1 & \text{for $s=0$,} \\
        x_d = \Bar{\bm{r}} \langle X \rangle & \text{for $s=1$.}
    \end{cases}
\end{equation}
In this way, it is ensured that one of the two decomposed branches cannot satisfy a selected $\Bar{\bm{r}}$.

For instance, consider the aforementioned example with $\Bar{\mathcal{R}}_1=\{\gamma_1, \Bar{\gamma}_3\}$, where $\gamma_1=\left((x_1^1,x_2^1), 1\right)$ and $\Bar{\gamma}_3=\left((x_2^0), 0\right)$. In the set, $X_2$ appears twice, $X_1$ once, and $X_3$ never. Thus, the component events are ordered as $X_2$, $X_1$, and $X_3$. Then, the conditional probabilities of the reduced rules are calculated as $P(\gamma_1 | b_1) \propto 0.9 \cdot 0.8 = 0.72$ and $P(\Bar{\gamma}_3 | b_1) \propto 0.2$. This orders the rules as $\gamma_1$ and $\Bar{\gamma_3}$. Then, one starts from $X_2$ and checks $\bm{r}_1 \langle X_2 \rangle = 1$, which satisfies $\bm{l}_1\langle X \rangle = 0 < \bm{r}_1 \langle X \rangle \leq \bm{u}_1\langle X \rangle = 1$. Therefore, $x_2^1$ is selected to decompose $b_1$. 

\subsubsection{Decomposition and new branches}

Once $x_d$ is decided, a branch $b=(\bm{l}, \bm{u}, s_d, s_u, p)$ can be decomposed into a lower branch $b_l=(\bm{l}, \bm{u}', s_d, s_u', p)$ and an upper branch $b_u=(\bm{l}', \bm{u}, s_d', s_u, p)$, where for $X\in \bm{X}$,
\begin{subequations}
\begin{equation}
    \bm{u}'\langle X \rangle = 
    \begin{cases}
        \bm{u}\langle X \rangle & \text{if $X \neq X_d$, } \\
        x_d - 1 & \text{if $X = X_d$}.
    \end{cases}
\end{equation}
\begin{equation}
    \bm{l}'\langle X \rangle = 
    \begin{cases}
        \bm{l}\langle X \rangle & \text{if $X \neq X_d$, } \\
        x_d & \text{if $X = X_d$}.
    \end{cases}
\end{equation}
\end{subequations}

Once new branches $b_l$ and $b_u$ are created, their $s_l$ and $s_u$ need to be decided by comparing the new bounds (i.e. $\bm{l}$ and $\bm{u}$) with a current rule set $\mathcal{R}$. The system state corresponding to a bound $\bm{x}$ is decided as $s=0$ if $\mathcal{R}_f$ has a failure rule that satisfies
\begin{equation}
    \bm{r}\langle X \rangle \geq \bm{x}\langle X \rangle, \ \forall X \in Scope[\gamma],
\end{equation}
and $s=1$ if $\mathcal{R}_s$ has a rule satisfying
\begin{equation}
    \bm{r}\langle X \rangle \leq \bm{x}\langle X \rangle, \ \forall X \in Scope[\gamma].
\end{equation}
When there are no rules that satisfy either of the conditions above, the system state is set to $s=\text{u}$. Meanwhile, the probabilities of new branches can be calculated by \eqref{eq:braPro} or approximated by \eqref{eq:braProApr}. 

For instance, in the example of Figure~\ref{subfig:toy}, system states corresponding to $\bm{l}_1 = (0,0,0)$ and $\bm{u}_1 = (1,1,0)$ are inferred as $s_{l,1}=0$ and $s_{u,1}=1$ from $\gamma_3$ and $\gamma_1$, respectively. The probability of the branch is calculated as $P(X_3\leq0)=0.3$.

\subsubsection{Termination}

When a for-loop is completed, the decomposition algorithm may terminate for two reasons. First, it terminates when a current branch set $\mathcal{B}_{\text{new}}$ contains only specified branches, which indicates that now the exact value of a system failure probability can be computed. Second, it terminates when all unspecified branches in $\mathcal{B}_{\text{new}}$ have an empty $\Bar{\mathcal{R}}$. This means that any further decomposition cannot identify more specified rules, and thus the rule evaluation procedure needs to be run to find a new rule. If neither of these conditions is met, another for-loop is performed for further decomposition.

\subsection{Rule evaluation}\label{sec:rule_eval}

\subsubsection{Inputs/outputs, strategies, and initialisation}

The rule evaluation module, illustrated on the upper left-hand side of Fig.~\ref{fig:brc}, takes as inputs a system performance function $\Phi(\cdot)$ and, from the decomposition module, a current set of branches, $\mathcal{B}$. It returns an updated rule set $\mathcal{R}$, which is fed back to the decomposition module.

The process is designed to select a next component vector state for system evaluation such that a new rule can be found. For the initialisation of the entire algorithm, a current rule set is defined as $\mathcal{R}=\emptyset$.

\subsubsection{Selecting a component vector state for system evaluation}\label{sec:x_star}

Given a branch set $\mathcal{B}$, we propose the following strategy for selecting a component vector state $\bm{x}^*$ with which to run a system performance function evaluation. First, the branches in $\mathcal{B}$ are sorted in descending order of their probabilities. Then, the algorithms selects the upper bound of the first branch having $s_u = \mathrm{u}$. If all branches have $s_u\neq\text{u}$, the algorithm then selects the lower bound of the first branch having $s_l = \mathrm{u}$. Through extensive numerical experiments, we found that considering upper bounds first is most effective in identifying rules with higher probabilities. This is because survival states in general have a higher probability than failure states, which makes it more likely for upper bounds to be associated with high-probability rules than lower bounds.

By selecting a $\bm{x}^*$ whose system state is unknown with the current rule set $\mathcal{R}$, it is guaranteed that the rule associated to $\bm{x}^*$ is a new addition to $\mathcal{R}$.

\subsubsection{Obtaining a new rule}\label{sec:sysFunUseInp}

The evaluation of $\Phi(\bm{x}^*)$ leads to a new rule. However, it can be more efficient to instead derive a minimal or sub-minimal (i.e. close to but not minimal) rule, which might be obtained as a by-product of a system simulation result. When this optional output is provided, the efficiency of the BRC algorithm is significantly improved. 

Optional outputs are not always available. In these cases, an analysed component vector state can be used as an associated rule although it is less likely to be (sub-)minimal. For instance, running a connectivity analysis on $(x_1^1, x_2^1,x_3^0)$ and confirming system survival, one can identify an associated survival rule $\left( (x_1^1, x_2^1), 1 \right)$, where the state $x_3^0$ is excluded since it is clear that this worst state does not contribute to system survival. Another example is that analysing $(x_1^0, x_2^0, x_3^1)$ and confirming system failure, one can infer a failure rule $\left( (x_1^0,x_2^0), 0 \right)$, where the best state $x_3^1$ is excluded for the same reason. We note that the rule is not minimal (as $x_2^0$ can be excluded to decide system failure) and therefore may be excluded later for analysis when a rule $\left( (x_1^0), 0 \right)$ is found (c.f. Section~\ref{sec:upd_rul}).

As an example of finding minimal rules, consider a network connectivity problem. In this case, a shortest path algorithm can be employed to find a shortest path between an origin and a destination node. It concludes system survival if a path is found and failure otherwise. Then, when a system state is survival, the obtained shortest path implies a minimal survival rule \cite{LiHe02,LimSon12,LeeSon21}. For instance, consider the example in Fig.~\ref{subfig:toy} and a branch $b_2=\left( (1,0,0), (1,1,0), \text{u}, \text{u}, 0.27 \right)$. When a shortest path algorithm is run on the upper bound, i.e. $(x_1^1,x_2^1,x_3^0)$, it would find a path $e_1 - e_2$ (as $e_3$ is unavailable). From this result, one can infer a minimal survival rule $\left( (x_1^1,x_2^1), 1 \right)$. Similarly, failure rules can be identified by employing a maximum flow algorithm, from which a minimum cut can be found \cite{AhuMagOrl93}. For instance, suppose that a maximum flow algorithm is run on the example lower bound, i.e. $(x_1^1,x_2^0,x_3^0)$. The algorithm will return a maximum flow as 0, concluding a system failure event, and a minimum cut $(e_2,e_3)$. This result can be used to infer a corresponding minimal failure rule $\left( (x_2^0,x_3^0), 0 \right)$.

A component vector state may be associated with more than one minimal rule. For example, in the system of Fig.~\ref{subfig:toy}, a component vector $(x_1^1,x_2^1,x_3^1)$ satisfies both survival rules $\left( (x_1^1,x_2^1), 1 \right)$ and $\left( (x_1^1,x_3^1), 1 \right)$. In case of an early termination of a branch-and-bound algorithm, narrowing down a bound on system failure probability can be expedited if rules with higher probabilities can be identified first \cite{LimSon12,ByuStr23}, where such a probability is calculated for a rule $\gamma=(\bm{r}, s)$ as
\begin{equation}
    P(\gamma) = 
    \begin{cases}
    P(\cap_{X\in Scope[\gamma]} \{ X \leq \bm{r}\langle X \rangle \} ) & \text{for $s=0$}, \\
        P(\cap_{X\in Scope[\gamma]} \{ X \geq \bm{r}\langle X \rangle \} ) & \text{for $s=1$}.
    \end{cases}
\end{equation}
A system function can be designed to identify rules with higher probabilities first.

The most well explored example is network connectivity, for which the selective Recursive Decomposition Algorithm (sRDA) has been proposed \cite{LimSon12}. The sRDA applies a shortest path algorithm by modifying the length of each edge to $-\log\left(P(x_n^1)\right)$, where $X_n$, $n=1,\ldots,N$, is a random variable standing for each edge \cite{LimSon12}. In this way, an obtained survival rule $\gamma$, being inferred from an obtained shortest path, will have the lowest value of $\sum_{X_n\in Scope[\gamma]} - \log \left( P(x^1_n) \right)$ and thus the greatest value of $\prod_{X_n\in Scope[\gamma]} P(x^1_n)$.

There are no universal strategies for identifying (sub-)minimal rules from a result of a system performance function evaluation. Such strategies need to be designed for a given system performance function. Multiple algorithms have been proposed to this end, including \cite{LiHe02,LimSon12,LeeSon21} for two-terminal network connectivity and \cite{JanLai08} for two-terminal maximum flow. We note that identifying associated survival rules is mostly straightforward, while finding failure rules may be more challenging. This is because in case of system survival, system analysis results often provide involved components in a system's operation, which can be considered components of a survival rule. In contrast, components responsible for a system failure are often less apparent. In the numerical examples of Section~\ref{sec:ex}, we assume that identifying (sub-)minimal survival rules is possible while finding (sub-)minimal failure rules is not.

\subsubsection{Updating a set of rules}\label{sec:upd_rul}

When a new rule $\gamma^*$ is obtained, the BRC algorithm updates a current $\mathcal{R}$ to add $\gamma^*$ and remove rules dominated by $\gamma^*$. A rule $\gamma_1$ is {\it dominated} by a rule $\gamma_2$ if $\gamma_1$ is always satisfied when $\gamma_2$ is satisfied (i.e. $\gamma_1$ represents component vector states that are a subset of $\gamma_2$). Formally, a rule $\gamma_1=(\bm{r}_1, s)$ is dominated by $\gamma_2=(\bm{r}_2, s)$ if $Scope[\gamma_2] \subset Scope[\gamma_1]$ and
\begin{subequations}
    \begin{equation}
         \bm{r}_2 \langle X \rangle \geq \bm{r}_1 \langle X \rangle, \ \forall X\in Scope[\gamma_2], \text{\quad for $s=0$, and}
    \end{equation}
    \begin{equation}
        \bm{r}_2 \langle X \rangle \leq \bm{r}_1 \langle X \rangle, \ \forall X\in Scope[\gamma_2], \text{\quad for $s=1$}.
    \end{equation}
\end{subequations}

\subsection{System failure probability evaluation}\label{sec:stop}

\subsubsection{Inputs and outputs}

Once a termination criterion is met in the decomposition algorithm (c.f. Section~\ref{sec:dec}), the BRC algorithm is terminated by evaluating the system failure probability in the reliability evaluation module, shown on the lower left-hand side of Fig.~\ref{fig:brc}. It takes as inputs the probability distribution of components $P(\bm{X})$ and a set of branches $\mathcal{B}$. Depending on a system scale, it returns either an exact value of system failure probability $p_f$, a bound $[p_{f,l},p_{f,u}]$, or an estimate $\hat{p}_f$ with a variance measure $\hat{\sigma}_f$. It also returns a rule set $\mathcal{R}$ and a branch set $\mathcal{B}$ as final outcomes of the BRC algorithm as they provide useful information on a given system event.

\subsubsection{Evaluation}
To obtain scalability to large-scale systems, the algorithm adopts one of three approaches to evaluating the system failure probability. The first approach is followed when $\mathcal{B}$ includes only specified branches. In this case, an exact value is readily obtained as
\begin{equation}
    p_f = \sum_{b\in\mathcal{B}_f} P(b),
\end{equation}
where $P(b)$ is calculated by \eqref{eq:braPro}. 

The second approach is followed is when a narrow enough bound is available, i.e. $P(\mathcal{B}_{\text{new},u}) \leq \epsilon_u$. In this case, the BRC algorithm returns a bound calculated as
\begin{equation}
    P(s^0) \in \left[ P(\mathcal{B}_f), 1-P(\mathcal{B}_s) \right] = \left[ \sum_{b\in\mathcal{B}_f} P(b), 1 - \sum_{b\in\mathcal{B}_s} P(b) \right],
\end{equation}  

The third approach is chosen when too many branches exist to continue decomposition, i.e. $\lvert \mathcal{B}_{\text{new}} \rvert \geq N_{b,\text{max}}$. In this case, we propose performing a hybrid inference by running Monte Carlo Simulation (MCS) over the set of unspecified branches $\mathcal{B}_u = \mathcal{B} \setminus (\mathcal{B}_f \cup \mathcal{B}_s)$. 
We utilise a Bayesian approach to quantify the uncertainty associated with this estimate. 
The prior of the system failure probability conditional on $\mathcal{B}_u$ is $\mathrm{Beta}(\alpha, \beta)$, i.e. the Beta distribution with parameters $\alpha$ and $\beta$ and bounds $P(\mathcal{B}_f)$ and $1-P(\mathcal{B}_s)$. $\alpha$ and $\beta$ are chosen to reflect prior knowledge. As a default, with no prior knowledge, \cite{BetPapStr22} shows that $\alpha=\beta=1$ is a conservative choice. Then, if one generates $M$ samples and observe $M_f$ failure samples, the posterior distribution follows $\mathrm{Beta}(\alpha + M_f, \beta + M - M_f)$, which has mean
\begin{equation}
    \hat{p}_{f,u} = (\alpha + M_f) / (\alpha + \beta + M)
\end{equation}
variance
\begin{equation}
    \hat{\sigma}^2_{f,u} = (\alpha + M_f) (\beta + M - M_f) / (\alpha + \beta + M)^2 (\alpha + \beta + M + 1)
\end{equation}
Combined with specified branches, the final estimate of the system failure probability is Beta-distributed with lower limit $P(\mathcal{B}_f)$, upper limit $P(\mathcal{B}_s)$, mean
\begin{equation}
    \hat{p}_{f} = P(\mathcal{B}_f) + P(\mathcal{B}_u) \cdot \hat{p}_{f,u}
\end{equation}
and variance
\begin{equation}\label{eq:pf_var}
    \hat{\sigma}^2_{f} = P(\mathcal{B}_u)^2 \cdot \hat{\sigma}^2_{f,u}.
\end{equation}

\subsubsection{Probability update}

In many instances, the failure probability estimate needs to be updated when new information on $\bm{X}$ becomes available. If all branches are specified, or if the bounds are sufficiently narrow, updating the system failure probability is straightforward with the posterior distribution of $\bm{X}$, $P'(\bm{X})$. No additional system performance function evaluations are required. Such evaluations can also be avoided when using sampling. In this case, the MCS analysis becomes importance sampling (IS). Each sample $\bm{x}_m$, $m=1,\ldots,M$, is assigned a weight $w_m = P(\bm{x}_m) / P'(\bm{x}_m)$. Then, the effective number of samples $M'$ and that of failure samples $M_f'$ is calculated as
\begin{subequations}\label{eq:upd}
    \begin{equation}
        M' = \sum_{m=1}^M w_m,
    \end{equation}
    \begin{equation}
        M_f' = \sum_{m=1}^M w_m\cdot \mathbb{I}(\Phi(\bm{x}_m) = 0),
    \end{equation}
\end{subequations}
by which a posterior distribution is evaluated as $\text{Beta}(\alpha + M_f', \beta + M' - M_f')$. IS can be ineffective when a sampling distribution is not optimised. However, when the posterior of $\bm{X}$ does not deviate too much from the prior, the latter is often a reasonable IS distribution. Additionally, the sampling variance is reduced by the presence of $\mathcal{B}_f$ and $\mathcal{B}_s$ as implied in \eqref{eq:pf_var}. The approach can thus be effective to update the reliability estimates. We support this argument by the numerical investigations in Section~\ref{sec:ema}.

\subsection{Didactic example}\label{sec:toy}

This section illustrates each step of the BRC algorithm using the example in Fig.~\ref{subfig:toy}. As before, we assume component failure probabilities $P(x_1^0)=0.1$, $P(x_2^0)=0.2$, and $P(x_3^0)=0.3$. The system failure event is defined as the disconnection of nodes 1 and 3. Following the discussion in Section~\ref{sec:sysFunUseInp}, we utilise connectivity analysis result to infer an associated minimal survival rule, while we assume that minimal failure rules are unavailable. The termination criteria set as $\epsilon_u=0$ and $N_{b,\text{max}}=\infty$, i.e. we aim for exact evaluation. The following iterations are summarised in Fig.~\ref{fig:toy_brc}. 

\subsubsection{Iteration 1.}
The algorithm begins with $\mathcal{R}=\emptyset$ and enters the decomposition algorithm. The decomposition enters the for-loop with $\mathcal{B}=\{b_0\}$ where $b_0=b_{\Omega}=\left((0,0,0),(1,1,1),\text{u},\text{u},1\right)$ and $\mathcal{B}_{\text{new}}=\emptyset$. No early termination criteria are met and $\mathcal{R}_0=\emptyset$ since $\mathcal{R}=\emptyset$. Therefore, $b_0$ is added to $\mathcal{B}_{\text{new}}$. Since this is the last element, the algorithm exits the for-loop. Since all unspecified rules in $\mathcal{B}_{\text{new}}$ (i.e. $b_0$) have an empty $\Bar{\mathcal{R}}$, the algorithm moves to the rule evaluation. It selects $\bm{u}_0$ for a next system simulation and learn that $s_{u,0}=1$ and this results in the survival rule $\left( (x_1^1, x_2^1), 1 \right)$. The current rule set is updated to $\mathcal{R} = \left\{ \left( (x_1^1, x_2^1), 1 \right) \right\} = \{ \gamma_1 \}$. 

\subsubsection{Iteration 2.}
The decomposition algorithm is entered with $\mathcal{R}=\{\gamma_1\}$. It enters the for-loop with $\mathcal{B}=\{b_0\}$ and $\mathcal{B}_{\text{new}}=\emptyset$. Starting with $b_0$, it proceeds to obtain $\Bar{\mathcal{R}}_0=\{\gamma_1\}$. $\Bar{\mathcal{R}}_0$ is non-empty and therefore $b_0$ can be decomposed. The components are ordered by their appearance in $\Bar{\mathcal{R}}_0$ as $X_1$ and $X_2$ (once) and $X_3$ (never). Then, since $\bm{r}_1\langle X_1 \rangle = 1$ can be used to decompose $b_0$, $x_1^1$ is selected for decomposition. This leads to new branches $b_1 = \left((0,0,0),(0,1,1),\text{u},\text{u},0.1\right)$ and $b_2 = \left((1,0,0),(1,1,1),\text{u},1,0.9\right)$. $s_{u,2}=1$ can be concluded by $\gamma_1$, while other system states remain unknown. The probabilities are calculated as $p_1 = P(X_1\leq0)=0.1$ and $p_2=P(X_1 > 0)=0.9$. These new branches are added to the new branch set, i.e. $\mathcal{B}_{\text{new}}=\{b_1,b_2\}$. There is no branch left in $\mathcal{B}$ and no termination criteria are met; the algorithm returns to the beginning of the for-loop.

The algorithm enters the for-loop with $\mathcal{B}=\{b_2,b_1\}$ and $\mathcal{B}_{\text{new}}=\emptyset$, where $\mathcal{B}$ is sorted by probabilities. Starting with $b_2$, the algorithm obtains $\Bar{\mathcal{R}}_2=\{\Bar{\gamma}_1\}$, where $\gamma_1$ is reduced to $\Bar{\gamma}_1 = \left( (x_2^1), 1 \right)$ since its condition $x_1^1$ is already satisfied by the lower bound $\bm{l}_2 = (1,0,0)$. Since $\Bar{\mathcal{R}}_2$ is non-empty, $b_2$ needs to be decomposed and the only remaining condition $x_2^1$ is selected. This splits $b_2$ into $b_3 = \left((1,0,0),(1,0,1),\text{u},\text{u},0.18\right)$ and $b_4 = \left((1,1,0),(1,1,1),1,1,0.72\right)$. Now the new branch set becomes $\mathcal{B}_{\text{new}} = \{ b_3, b_4 \}$. Both system states of $b_4$ can be confirmed as 1 by $\gamma_1$ as 1, i.e. $b_4$ is a specified branch. Moving to $b_1$, it has $\Bar{\mathcal{R}}_1=\emptyset$ as the condition $x_1^1$ of $\gamma_1$ cannot be satisfied by its upper bound $\bm{u}_1$. The new branch set is updated to $\mathcal{B}_{\text{new}} = \{ b_3, b_4, b_1 \}$. The algorithm reaches the end of the for-loop.

The algorithm returns to the beginning of the for-loop with the sorted branches $\mathcal{B} = \{ b_4, b_3, b_1 \}$. The specified branch $b_4$ is immediately added to $\mathcal{B}_{\text{new}}$. $b_3$ has an empty $\Bar{\mathcal{R}}_3$ as $x_2^1$ of $\gamma_1$ cannot be satisfied by its upper bound $\bm{u}_3$; it is added to $\mathcal{B}_{\text{new}}$. As previously evaluated, $b_1$ has $\Bar{\mathcal{R}}_1=\emptyset$. The for-loop is completed with $\mathcal{B}_{\text{new}} = \{ b_4, b_3, b_1 \}$. Since all unspecified branches in $\mathcal{B}_{\text{new}}$ have an empty $\Bar{\mathcal{R}}$, the algorithm goes back to the rule evaluation with $\mathcal{B} = \{ b_4, b_3, b_1 \}$. To run the system function, $\bm{u}_3$ is selected, which is the first upper bound with an unknown system state. The system function returns $s_{u,3}=1$ and a new rule $\gamma_2 = \left( (x_1^1, x_3^1), 1 \right)$. This updates the rule set to $\mathcal{R} = \left\{ \left( (x_1^1, x_2^1), 1 \right), \left( (x_1^1, x_3^1), 1 \right) \right\} = \{ \gamma_1, \gamma_2 \}$. 

\subsubsection{Iteration 3.}
The algorithm re-enters the for-loop with the reduced rules of $b_2$, $\Bar{\mathcal{R}}_2 = \{\Bar{\gamma}_1,\Bar{\gamma}_2\}$, where $\Bar{\gamma}_1=\left( (x_2^1), 1 \right)$ and $\Bar{\gamma}_2=\left( (x_3^1), 1 \right)$. For decomposition, the components are ordered as $X_2$ and $X_3$ (appearing once) and $X_1$ (never). Then, the rules are ordered as $\Bar{\gamma}_1$ and $\Bar{\gamma}_2$ because of their probabilities $P(\Bar{\gamma}_1 | b_2) \propto P(X_2 \geq 1) = 0.8$ and $P(\Bar{\gamma}_2 | b_2) \propto P(X_3 \geq 1) = 0.7$. Starting from $X_2$ and $\Bar{\gamma}_1$, $x_2^1$ is selected for decomposition, splitting $b_2$ into the same $b_3$ and $b_4$ as in Iteration 2. Moving to $b_1$, the branch has an empty $\Bar{\mathcal{R}}_1$ since the condition $x_1^1$ required by both $\gamma_1$ and $\gamma_2$ cannot be satisfied by $\bm{u}_1$. The for-loop is ended with $\mathcal{B}_{\text{new}} = \{ b_3, b_4, b_1 \}$.

Next, the algorithm re-enters the for-loop with the sorted branches $\mathcal{B} = \{ b_4, b_3, b_1 \}$. The specified branch $b_4$ is added to $\mathcal{B}_{\text{new}}$ and $b_3$ is analysed next. This leads to $\Bar{\mathcal{R}}_3 = \{ \Bar{\gamma_2} \}$ since $x_2^1$ required by $\gamma_1$ cannot be satisfied. The only condition state in $\Bar{\mathcal{R}}_3$, $x_3^1$ is selected for decomposition. This splits $b_3$ into $b_5=\left( (1,0,0), (1,0,0), \mathrm{u}, \mathrm{u}, 0.054 \right)$ and $b_6=\left( (1,0,1), (1,0,1), 1, 1, 0.126 \right)$, where $b_6$ has both system states known as 1 by $\gamma_2$. $b_1$ has an empty $\Bar{\mathcal{R}}_1$ and therefore is not decomposed. The new branch set becomes $\mathcal{B}_{\text{new}}=\{b_5, b_1 \}$, where all unspecified branches have an empty $\Bar{\mathcal{R}}$. The algorithm goes back to the rule evaluation.

The branches are sorted by their probabilities as $\mathcal{B}=\{b_1, b_5 \}$. The upper bound $\bm{l}_1$ is selected for a system function evaluation as the upper bound of the branch with the highest probability. The system function returns $s_{u,1}=0$ and infer a failure rule $\left((x_1^0), 0\right)$ where the states $x_2^1$ and $x_3^2$ are excluded (c.f. Section~\ref{sec:sysFunUseInp}). This updates the rule set to $\mathcal{R} = \left\{ \left( (x_1^1, x_2^1), 1 \right), \left( (x_1^1, x_3^1), 1 \right), \left((x_1^0), 0\right) \right\} = \{ \gamma_1, \gamma_2, \gamma_3 \}$. 

\subsubsection{Iteration 4.}
By the same procedure as in Iteration 3, the algorithm yields a branch set $\mathcal{B}=\{b_1,b_4,b_5,b_6\}$. $b_1$ is now a specified set since $s_{l,1} = s_{u,1} = 0$ can be inferred from $\gamma_3$. Therefore, $b_5$ is the only unspecified branch. $\bm{u}_5$ is selected to run a system performance function, returning $s_{u,5}=0$ and identify a new failure rule $\left( (x_2^0, x_3^0), 0\right)$. Since $\bm{u}_5=\bm{l}_5$, the algorithm can infer from the new failure rule that $s_{l,5}=0$, which specifies $b_5$ as failure branch. Since there is no remaining unspecified rule, the algorithm is terminated with $\mathcal{B}=\{b_1,b_4,b_5,b_6\}$ and $\mathcal{R}=\left\{ \left( (x_1^1, x_2^1), 1 \right), \left( (x_1^1, x_3^1), 1 \right), \left((x_1^0), 0\right), \left( (x_2^0, x_3^0), 0\right) \right\}$.

\subsubsection{Final inference.}
The system failure probability is calculated as $P(b_1)+P(b_5)=0.1+0.054=0.154$. The BRC algorithm runs the system performance function four times to identify all existing rules. 

\begin{figure}[h!]
    \centering
    \begin{subfigure}[b]{\textwidth}
        \centering
        \tikzstyle{block} = [draw, rectangle, minimum height=0.8cm, minimum width=1.4cm]

\begin{tikzpicture}[thick, main/.style = {draw, circle}, node distance=1.5cm] 

\newcommand\xdist{2}
\newcommand\ydist{1}
\newcommand\xtext{-1.4}
\newcommand\ytext{1.2}
\newcommand\twn{6cm}
\newcommand\ys{1cm} 

\node[main] at (0,0) (b0) {$b_0$};
\node[align=left] at (0,\ytext) {\textcolor{red}{$\bm{u}_0=(1,1,1)$, $s_{u,0}=\mathrm{u}$} \\ $\bm{l}_0=(0,0,0)$, $s_{l,0}=\mathrm{u}$ \\ $p_0=1$, $\Bar{\mathcal{R}}_0=\emptyset$};
\node[align=center] at (0,-\ys) {\textrightarrow \textcolor{red}{$s_{u,0}=1$}, $\mathcal{R}=\big\{\textcolor{red}{\big((x_1^1,x_2^1), 1\big)}\big\}$};
    
\end{tikzpicture}
        \caption{Iteration 1.}
        \label{fig:toy_s1}
    \end{subfigure}
    
    \medskip
    \begin{subfigure}[b]{\textwidth}
        \centering
        \tikzstyle{block} = [draw, rectangle, minimum height=0.8cm, minimum width=1.4cm]

\begin{tikzpicture}[thick, main/.style = {draw, circle}, node distance=1.5cm] 

\newcommand\xdist{1}
\newcommand\ydist{2.2}
\newcommand\xtext{2.7}
\newcommand\ytext{0.3}
\newcommand\ytexti{1.2}
\newcommand\ys{1cm} 
\newcommand\nsz{0.75cm}

\node[main, minimum size = \nsz] at (0,0) (b0) {$b_0$};
\node[align=left] at (0,\ytexti) {$\bm{u}_0=(1,1,1)$, $s_{u,0}=1$ \\ $\bm{l}_0=(0,0,0)$, $s_{l,0}=\text{u}$ \\ $p_0 = 1$, $\Bar{\mathcal{R}}_0=\{\gamma_1\}$};

\node[main] at (-\xdist,-\ydist) (b1) {$b_1$};
\draw[->] (b0) -- (b1) node[midway,above left] {$X_1 \leq 0$};
\node[align=left] at (-\xdist-\xtext, -\ydist+\ytext) {$\bm{u}_1=(0,1,1)$, $s_{u,1}=\mathrm{u}$ \\ $\bm{l}_1=(0,0,0)$, $s_{l,1}=\text{u}$ \\ $p_1 = 0.1$, $\Bar{\mathcal{R}}_1 = \emptyset$ };

\node[main] at (\xdist,-\ydist) (b2) {$b_2$};
\draw[->] (b0) -- (b2) node[midway,above right] {$X_1 \geq 1$};
\node[align=left] at (\xdist+\xtext, -\ydist+\ytext) {$\bm{u}_2=(1,1,1)$, $s_{u,2}=1$ \\ $\bm{l}_2=(1,0,0)$, $s_{l,2}=\mathrm{u}$ \\ $p_2 = 0.9$, $\Bar{\mathcal{R}}_2=\{\Bar{\gamma}_1\}$};

\node[main] at (-\xdist+\xdist,-\ydist-\ydist) (b3) {$b_3$};
\draw[->] (b2) -- (b3) node[midway,above left] {$X_2 \leq 0$};
\node[align=left] at (-\xdist+\xdist-\xtext, -\ydist-\ydist+\ytext) {\textcolor{red}{$\bm{u}_3=(1,0,1)$, $s_{u,3}=\mathrm{u}$} \\ $\bm{l}_3=(1,0,0)$, $s_{l,3}=\mathrm{u}$ \\ $p_3 = 0.18$, $\Bar{\mathcal{R}}_3=\emptyset$};

\node[main, fill=black!10] at (-\xdist+3*\xdist,-\ydist-\ydist) (b4) {$b_4$};
\draw[->] (b2) -- (b4) node[midway,above right] {$X_2 \geq 1$};
\node[align=left] at (-\xdist+3*\xdist+\xtext, -\ydist-\ydist+\ytext) {$\bm{u}_4=(1,1,1)$, $s_{u,4}=1$ \\ $\bm{l}_4=(1,1,0)$, $s_{l,4}=1$ \\ $p_4=0.72$};

\node[align=center] at (0,-\ydist-\ydist-5.3cm) {\textrightarrow \textcolor{red}{$s_{u,3}=1$}, $\mathcal{R}=\big\{\big((x_1^1,x_2^1), 1\big), \textcolor{red}{\big((x_1^1,x_3^1), 1\big)} \big\}$};

\end{tikzpicture}
        \caption{Iteration 2.}
        \label{fig:toy_s2}
    \end{subfigure}
\end{figure}
\begin{figure}\ContinuedFloat
    \centering
    \begin{subfigure}[b]{\textwidth}
        \centering
        \tikzstyle{block} = [draw, rectangle, minimum height=0.8cm, minimum width=1.4cm]

\begin{tikzpicture}[thick, main/.style = {draw, circle}, node distance=1.5cm] 

\newcommand\xdist{1}
\newcommand\ydist{2.2}
\newcommand\xtext{2.7}
\newcommand\ytext{0.3}
\newcommand\ytexti{1.2}
\newcommand\ys{1cm} 
\newcommand\nsz{0.75cm}

\node[main, minimum size = \nsz] at (0,0) (b0) {$b_0$};
\node[align=left] at (0,\ytexti) {$\bm{u}_0=(1,1,1)$, $s_{u,0}=1$ \\ $\bm{l}_0=(0,0,0)$, $s_{l,0}=\text{u}$ \\ $p_0 = 1$, $\Bar{\mathcal{R}}_0=\{\gamma_1, \gamma_2\}$};

\node[main] at (-\xdist,-\ydist) (b1) {$b_1$};
\draw[->] (b0) -- (b1) node[midway,above left] {$X_1 \leq 0$};
\node[align=left] at (-\xdist-\xtext, -\ydist+\ytext) {\textcolor{red}{ $\bm{u}_1=(0,1,1)$, $s_{u,1}=\mathrm{u}$} \\ $\bm{l}_1=(0,0,0)$, $s_{l,1}=\text{u}$ \\ $p_1 = 0.1$, $\Bar{\mathcal{R}}_1 = \emptyset$ };

\node[main] at (\xdist,-\ydist) (b2) {$b_2$};
\draw[->] (b0) -- (b2) node[midway,above right] {$X_1 \geq 1$};
\node[align=left] at (\xdist+\xtext, -\ydist+\ytext) {$\bm{u}_2=(1,1,1)$, $s_{u,2}=1$ \\ $\bm{l}_2=(1,0,0)$, $s_{l,2}=\mathrm{u}$ \\ $p_2 = 0.9$, $\Bar{\mathcal{R}}_2=\{\Bar{\gamma}_1, \Bar{\gamma}_2\}$ };

\node[main] at (-\xdist+\xdist,-\ydist-\ydist) (b3) {$b_3$};
\draw[->] (b2) -- (b3) node[midway,above left] {$X_2 \leq 0$};
\node[align=left] at (-\xdist+\xdist-\xtext, -\ydist-\ydist+\ytext) {$\bm{u}_3=(1,0,1)$, $s_{u,3}=1$ \\ $\bm{l}_3=(1,0,0)$, $s_{l,3}=\mathrm{u}$ \\ $p_3 = 0.18$,  $\Bar{\mathcal{R}}_3=\{\Bar{\gamma}_2\}$};

\node[main, fill=black!10] at (-\xdist+3*\xdist,-\ydist-\ydist) (b4) {$b_4$};
\draw[->] (b2) -- (b4) node[midway,above right] {$X_2 \geq 1$};
\node[align=left] at (-\xdist+3*\xdist+\xtext, -\ydist-\ydist+\ytext) {$\bm{u}_4=(1,1,1)$, $s_{u,4}=1$ \\ $\bm{l}_4=(1,1,0)$, $s_{l,4}=1$ \\ $p_4=0.72$};

\node[main] at (-\xdist,-\ydist-\ydist-\ydist) (b5) {$b_5$};
\draw[->] (b3) -- (b5) node[midway,above left] {$X_3 \leq 0$};
\node[align=left] at (-\xdist-\xtext, -\ydist-\ydist-\ydist+\ytext) {$\bm{u}_5=(1,0,0)$, $s_{u,5}=\text{u}$ \\ $\bm{l}_5=(1,0,0)$, $s_{l,5}=\mathrm{u}$ \\ $p_5 = 0.054$, $\Bar{\mathcal{R}}_5=\emptyset$};

\node[main, fill=black!10] at (-\xdist+2*\xdist,-\ydist-\ydist-\ydist) (b6) {$b_6$};
\draw[->] (b3) -- (b6) node[midway,above right] {$X_3 \geq 1$};
\node[align=left] at (-\xdist+2*\xdist+\xtext, -\ydist-\ydist-\ydist+\ytext) {$\bm{u}_6=(1,0,1)$, $s_{u,6}=1$ \\ $\bm{l}_6=(1,0,1)$, $s_{l,6}=1$ \\ $p_6=0.126$};

\node[align=center] at (0,-4\ydist-6.3cm) {\textrightarrow \textcolor{red}{$s_{u,1}=0$}, $\mathcal{R}=\big\{\big((x_1^1,x_2^1), 1\big), \big((x_1^1,x_3^1), 1\big), \textcolor{red}{\big((x_1^0), 0\big)} \big\}$};

\end{tikzpicture}
        \caption{Iteration 3.}
        \label{fig:toy_s3}
    \end{subfigure}
\end{figure}

\begin{figure}\ContinuedFloat
    \begin{subfigure}[b]{\textwidth}
        \centering
        \tikzstyle{block} = [draw, rectangle, minimum height=0.8cm, minimum width=1.4cm]

\begin{tikzpicture}[thick, main/.style = {draw, circle}, node distance=1.5cm] 

\newcommand\xdist{1}
\newcommand\ydist{2.2}
\newcommand\xtext{2.7}
\newcommand\ytext{0.3}
\newcommand\ytexti{1.2}
\newcommand\ys{1cm} 
\newcommand\nsz{0.75cm}

\node[main, minimum size = \nsz] at (0,0) (b0) {$b_0$};
\node[align=left] at (0,\ytexti) {$\bm{u}_0=(1,1,1)$, $s_{u,0}=1$ \\ $\bm{l}_0=(0,0,0)$, $s_{l,0}=\text{u}$ \\ $p_0 = 1$, $\Bar{\mathcal{R}}_0=\{\gamma_1, \gamma_2\}$};

\node[main, fill=black!10] at (-\xdist,-\ydist) (b1) {$b_1$};
\draw[->] (b0) -- (b1) node[midway,above left] {$X_1 \leq 0$};
\node[align=left] at (-\xdist-\xtext, -\ydist+\ytext) {$\bm{u}_1=(0,1,1)$, $s_{u,1}=0$ \\ $\bm{l}_1=(0,0,0)$, $s_{l,1}=0$ \\ $p_1 = 0.1$ };

\node[main] at (\xdist,-\ydist) (b2) {$b_2$};
\draw[->] (b0) -- (b2) node[midway,above right] {$X_1 \geq 1$};
\node[align=left] at (\xdist+\xtext, -\ydist+\ytext) {$\bm{u}_2=(1,1,1)$, $s_{u,2}=1$ \\ $\bm{l}_2=(1,0,0)$, $s_{l,2}=\mathrm{u}$ \\ $p_2 = 0.9$, $\Bar{\mathcal{R}}_2=\{\Bar{\gamma}_1,\Bar{\gamma}_2\}$ };

\node[main] at (-\xdist+\xdist,-\ydist-\ydist) (b3) {$b_3$};
\draw[->] (b2) -- (b3) node[midway,above left] {$X_2 \leq 0$};
\node[align=left] at (-\xdist+\xdist-\xtext, -\ydist-\ydist+\ytext) {$\bm{u}_3=(1,0,1)$, $s_{u,3}=1$ \\ $\bm{l}_3=(1,0,0)$, $s_{l,3}=\mathrm{u}$ \\ $p_3 = 0.18$,  $\Bar{\mathcal{R}}_3=\{\Bar{\gamma}_2\}$};

\node[main, fill=black!10] at (-\xdist+3*\xdist,-\ydist-\ydist) (b4) {$b_4$};
\draw[->] (b2) -- (b4) node[midway,above right] {$X_2 \geq 1$};
\node[align=left] at (-\xdist+3*\xdist+\xtext, -\ydist-\ydist+\ytext) {$\bm{u}_4=(1,1,1)$, $s_{u,4}=1$ \\ $\bm{l}_4=(1,1,0)$, $s_{l,4}=1$ \\ $p_4=0.72$};

\node[main] at (-\xdist,-\ydist-\ydist-\ydist) (b5) {$b_5$};
\draw[->] (b3) -- (b5) node[midway,above left] {$X_3 \leq 0$};
\node[align=left] at (-\xdist-\xtext, -\ydist-\ydist-\ydist+\ytext) {\textcolor{red}{ $\bm{u}_5=(1,0,0)$, $s_{u,5}=\text{u}$} \\ $\bm{l}_5=(1,0,0)$, $s_{l,5}=\mathrm{u}$ \\ $p_5 = 0.054$, $\Bar{\mathcal{R}}_5=\emptyset$};

\node[main, fill=black!10] at (-\xdist+2*\xdist,-\ydist-\ydist-\ydist) (b6) {$b_6$};
\draw[->] (b3) -- (b6) node[midway,above right] {$X_3 \geq 1$};
\node[align=left] at (-\xdist+2*\xdist+\xtext, -\ydist-\ydist-\ydist+\ytext) {$\bm{u}_6=(1,0,1)$, $s_{u,6}=1$ \\ $\bm{l}_6=(1,0,1)$, $s_{l,6}=1$ \\ $p_6=0.126$};

\node[align=center] at (0,-4\ydist-6.3cm) {\textrightarrow \textcolor{red}{$s_{u,5}=0$}, $\mathcal{R}=\big\{\big((x_1^1,x_2^1), 1\big), \big((x_1^1,x_3^1), 1\big), \big((x_1^0), 0\big), \textcolor{red}{\big((x_2^0, x_3^0), 0\big)} \big\}$};

\end{tikzpicture}
        \caption{Iteration 4 (last).}
        \label{fig:toy_s4}
    \end{subfigure}

    \caption{Procedures of the BRC algorithm applied for the system in Fig.~\ref{subfig:toy}. System failure probability is defined as disconnectivity between node 1 and 3. A component event $X_n$, $n=1,2,3$ represents the binary state of edge $e_n$, where $P(x_1^0)=0.1$, $P(x_2^0)=0.2$, and $P(x_3^0)=0.3$. In the figures, specified branches are marked gray and their $\Bar{\mathcal{R}}$ are not specified as they do not require further decomposition. At each iteration, the selected state for system simulation is marked red; the system simulation result and the newly obtained rule is also marked red.}
    \label{fig:toy_brc}
\end{figure}
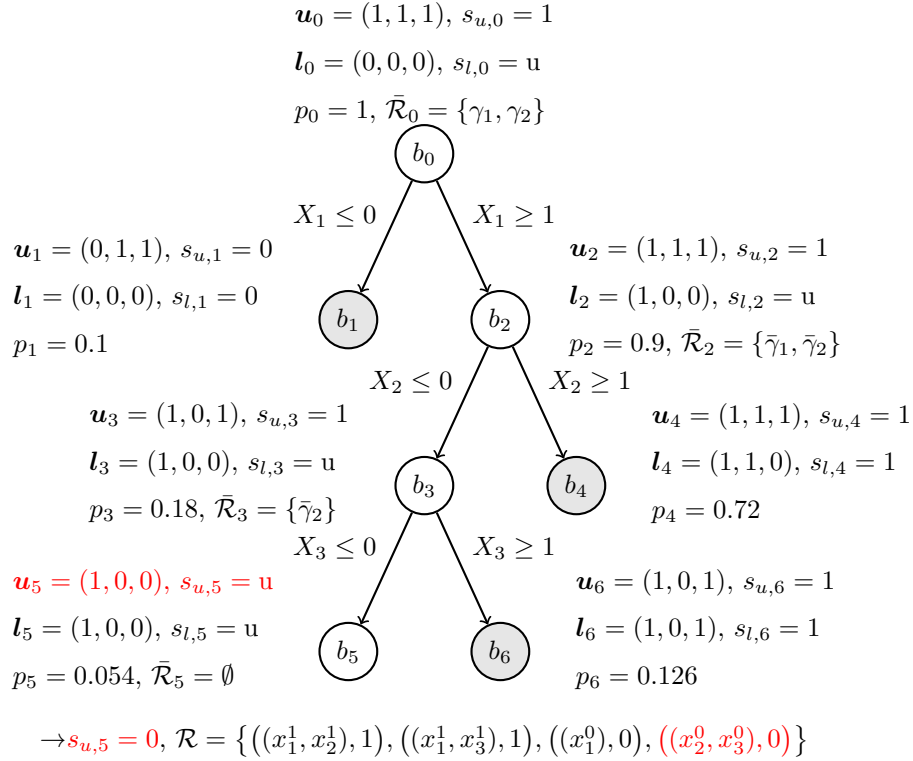
\FloatBarrier

\section{Numerical investigations}\label{sec:ex}

\subsection{Two-terminal maximum flow}\label{sec:mf}

We investigate a network whose performance is modeled by two-terminal maximum flow taken from \cite{JanLai08}.
This enables a comparison of the BRC algorithm to a state-of-the-art branch-and-bound algorithm, specifically the one of Jane and Laih \cite{JanLai08}. We note that the latter is applicable only to two-terminal maximum flow reliability, while the BRC algorithm can be used with any coherent system.  

\begin{figure}[h!]
    \centering
    \includegraphics[width=0.6\textwidth]{ 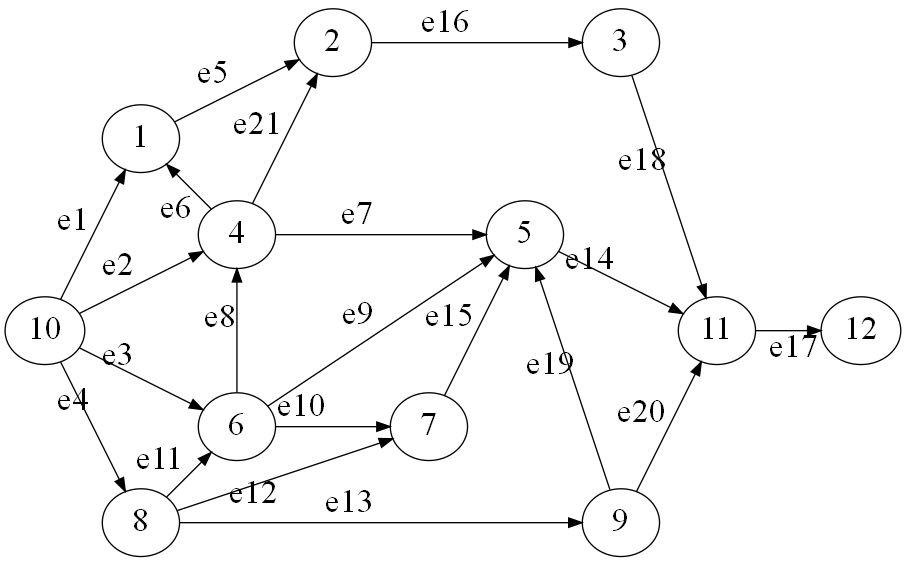}
    \caption{The network for two-terminal maximum flow reliability analysed in \cite{JanLai08}. The origin and the destination nodes are 10 and 12, respectively.}
    \label{fig:mf_net}
\end{figure}

Fig.~\ref{fig:mf_net} shows the example network analysed in \cite{JanLai08}, which consists of 12 nodes and 21 edges. The random variables $X_n$, $n=1,\ldots,21$, represent the state of a corresponding edge $e_n$. $X_n$ has three states 0, 1, and 2, which represent $e_n$ having a flow capacity of 0, 3, and 5, respectively. Their probability distributions $P(X_n)$, $n=1,\ldots,22$ are the same as in \cite{JanLai08}, where they are considered statistically independent. The system failure event is defined as the maximum flow from node 10 to 12 should be equal to or less than $d$. The system function is based on maximum flow analysis. In case of a system survival event, it returns the states of the components (i.e. edges) that constitute the flow path as a survival rule. For instance, given $d=3$, if a maximum flow is obtained with a path $e_3-e_9-e_{14}-e_{17}$ and their input states are $x_3^1$, $x_9^2$, $x_{14}^1$, and $x_{17}^2$, it returns a survival rule $\left( \left( x_3^1, x_9^2, x_{14}^1, x_{17}^2 \right), 1 \right)$. On the other hand, when a failure system state is obtained, the algorithm does not return any rule because identifying a (sub-)minimal failure rule is not straightforward. 

\begin{figure}[h!]
    \centering
    \begin{subfigure}[b]{0.48\textwidth}
    \centering
        \includegraphics[width=\textwidth]{ 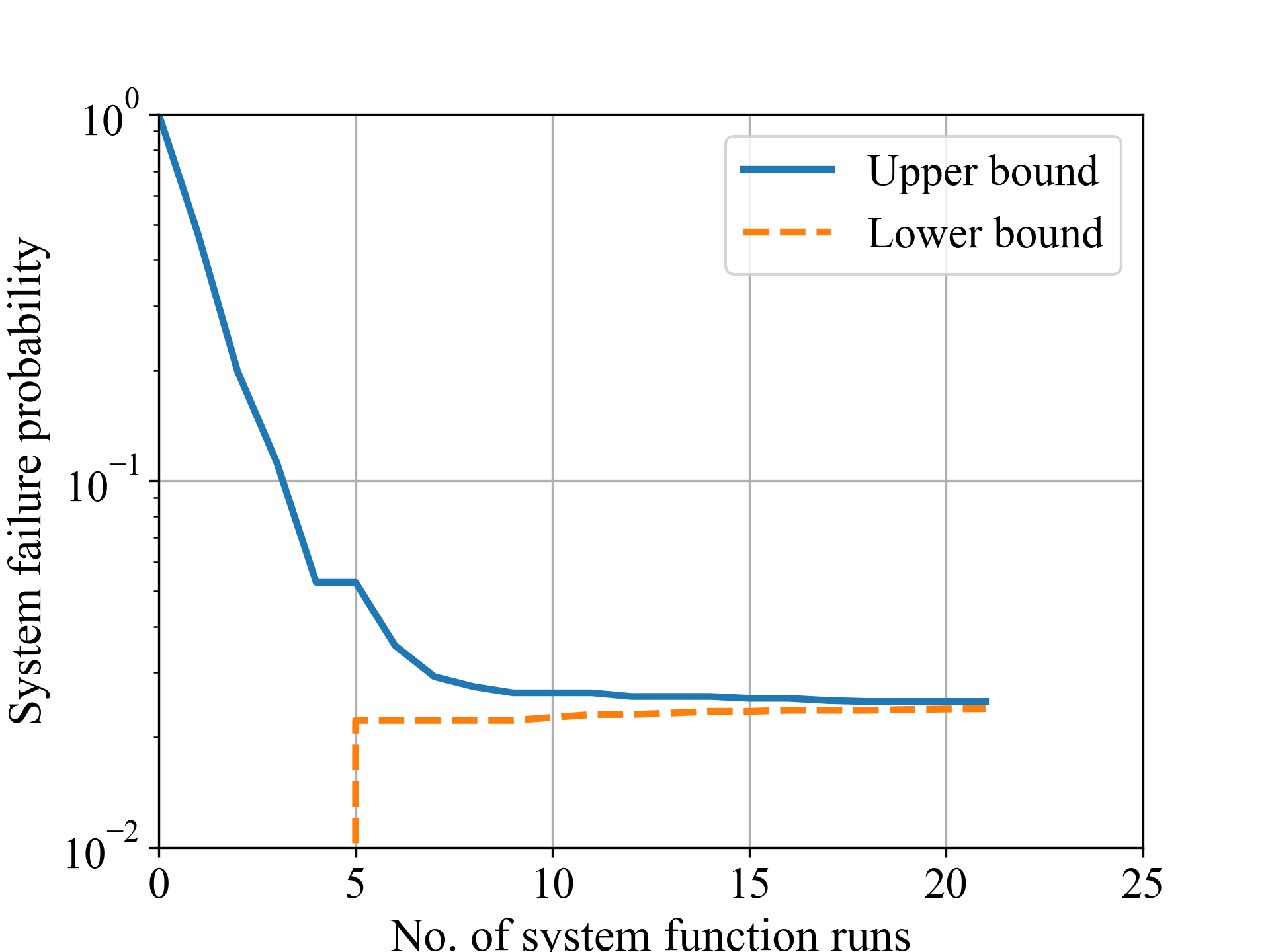}
        \caption{Convergence of failure probability (case $d=1,2,3$).}
        \label{fig:pf_d1}
    \end{subfigure}
    \hfill
    \begin{subfigure}[b]{0.48\textwidth}
    \centering
        \includegraphics[width=\textwidth]{ 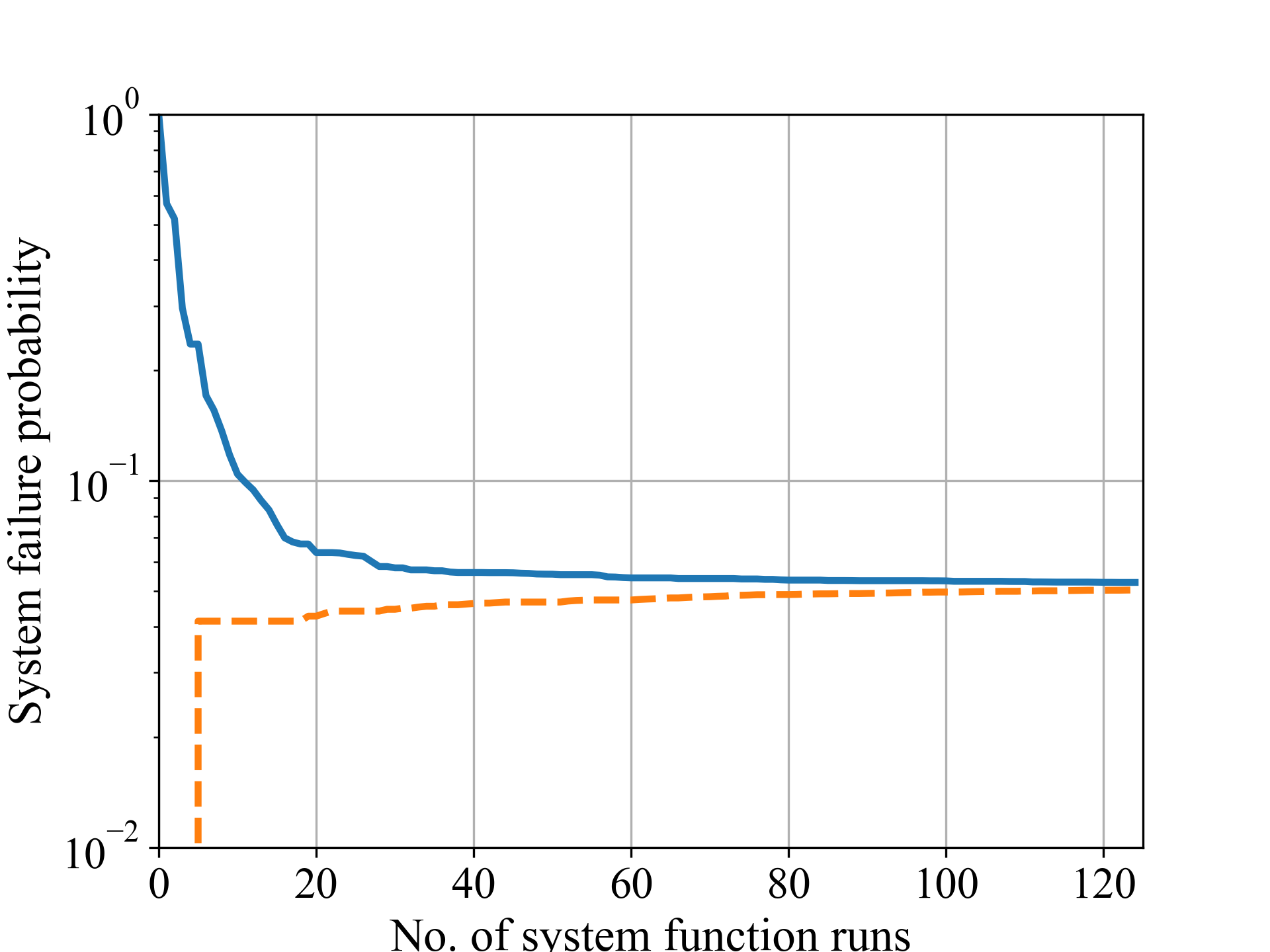}
        \caption{Convergence of failure probability (case $d=4,5$).}
        \label{fig:pf_d4}
    \end{subfigure}
    
    \centering
    \begin{subfigure}[b]{0.48\textwidth}
    \centering
        \includegraphics[width=\textwidth]{ 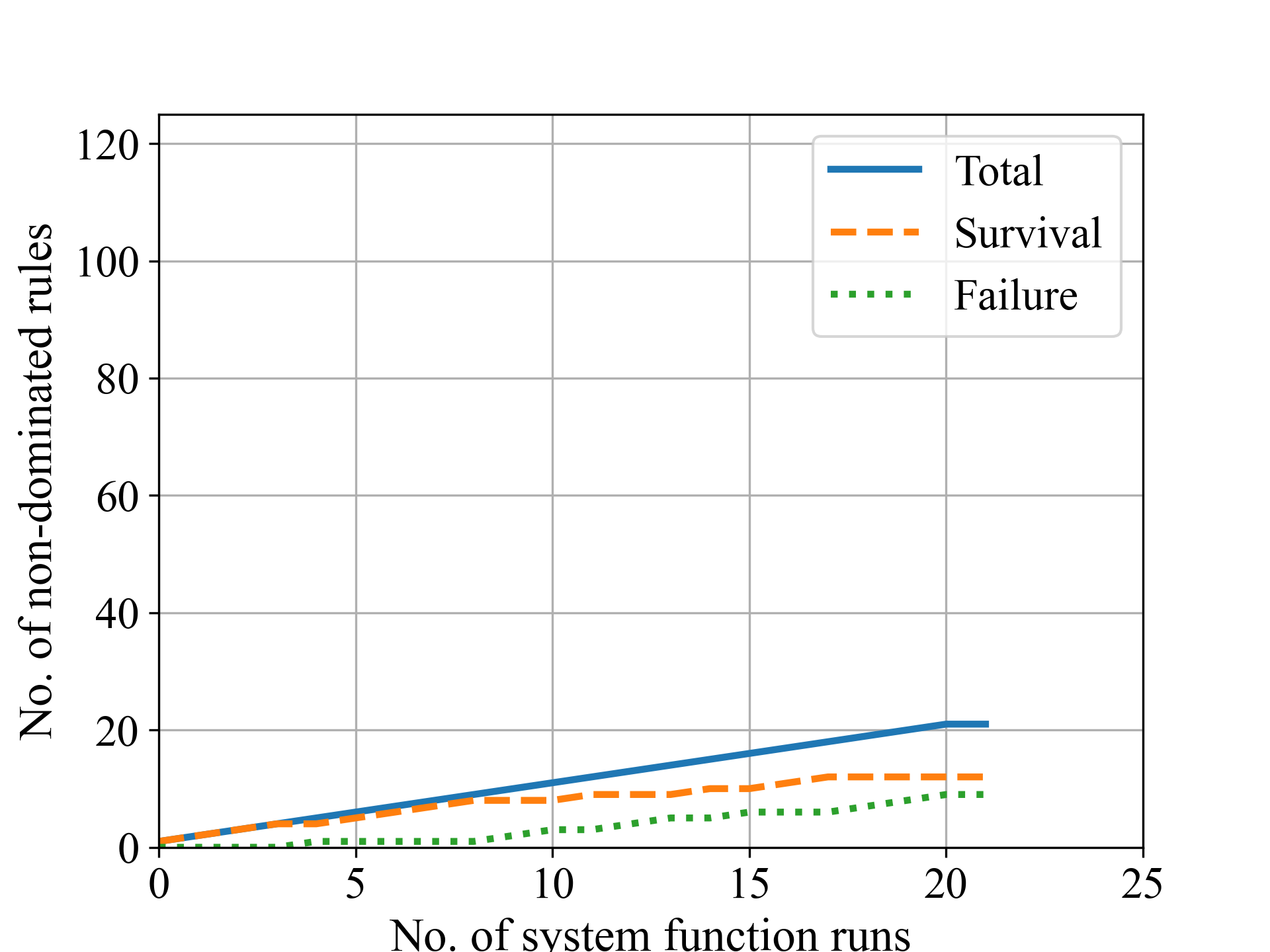}
        \caption{Number of non-dominated rules (case $d=1,2,3$).}
        \label{fig:ru_d1}
    \end{subfigure}
    \hfill
    \begin{subfigure}[b]{0.48\textwidth}
    \centering
        \includegraphics[width=\textwidth]{ 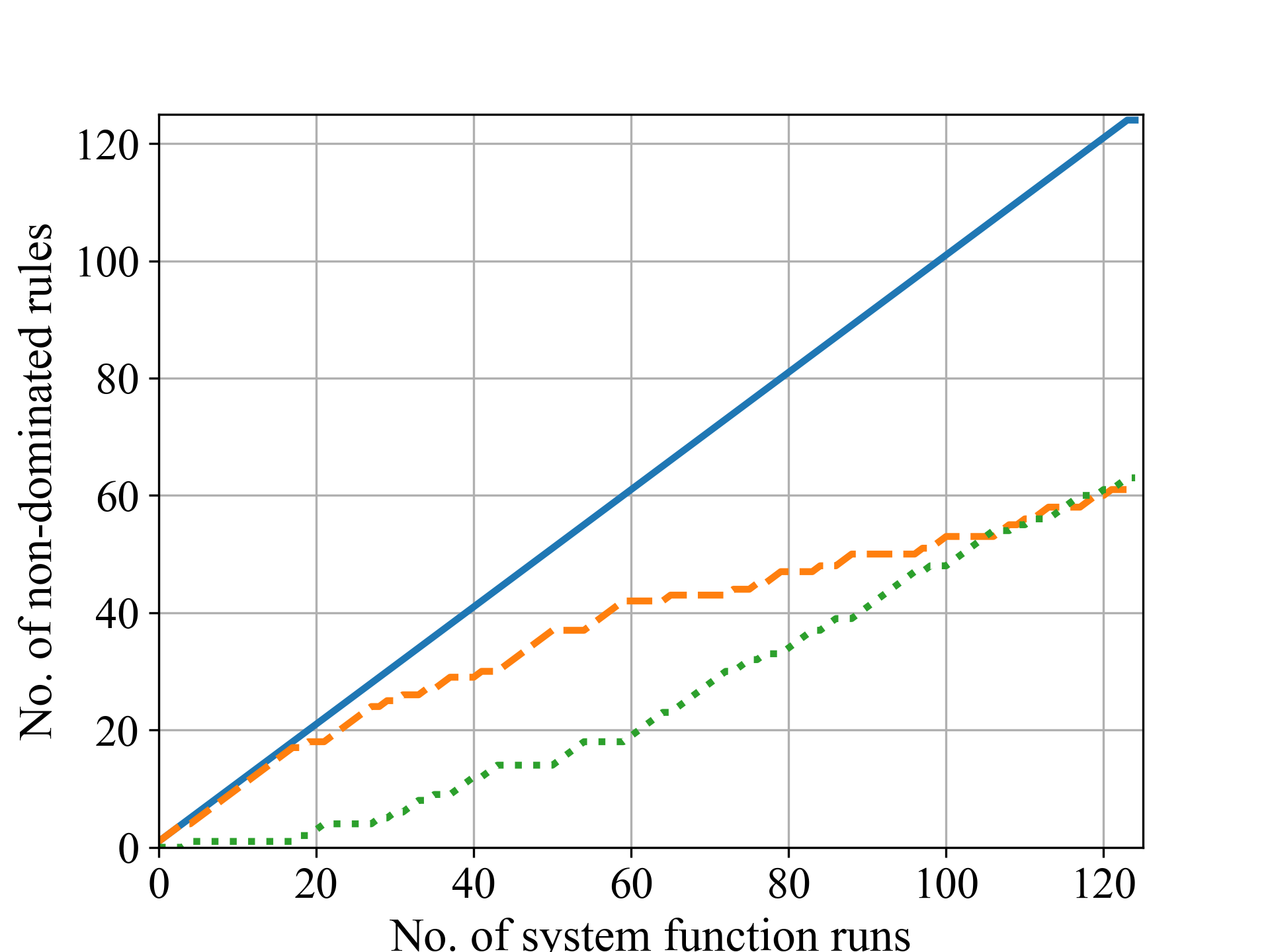}
        \caption{Number of non-dominated rules (case $d=4,5$).}
        \label{fig:ru_d4}
    \end{subfigure}
    
    \begin{subfigure}[b]{0.48\textwidth}
        \centering
        \includegraphics[width=\textwidth]{ 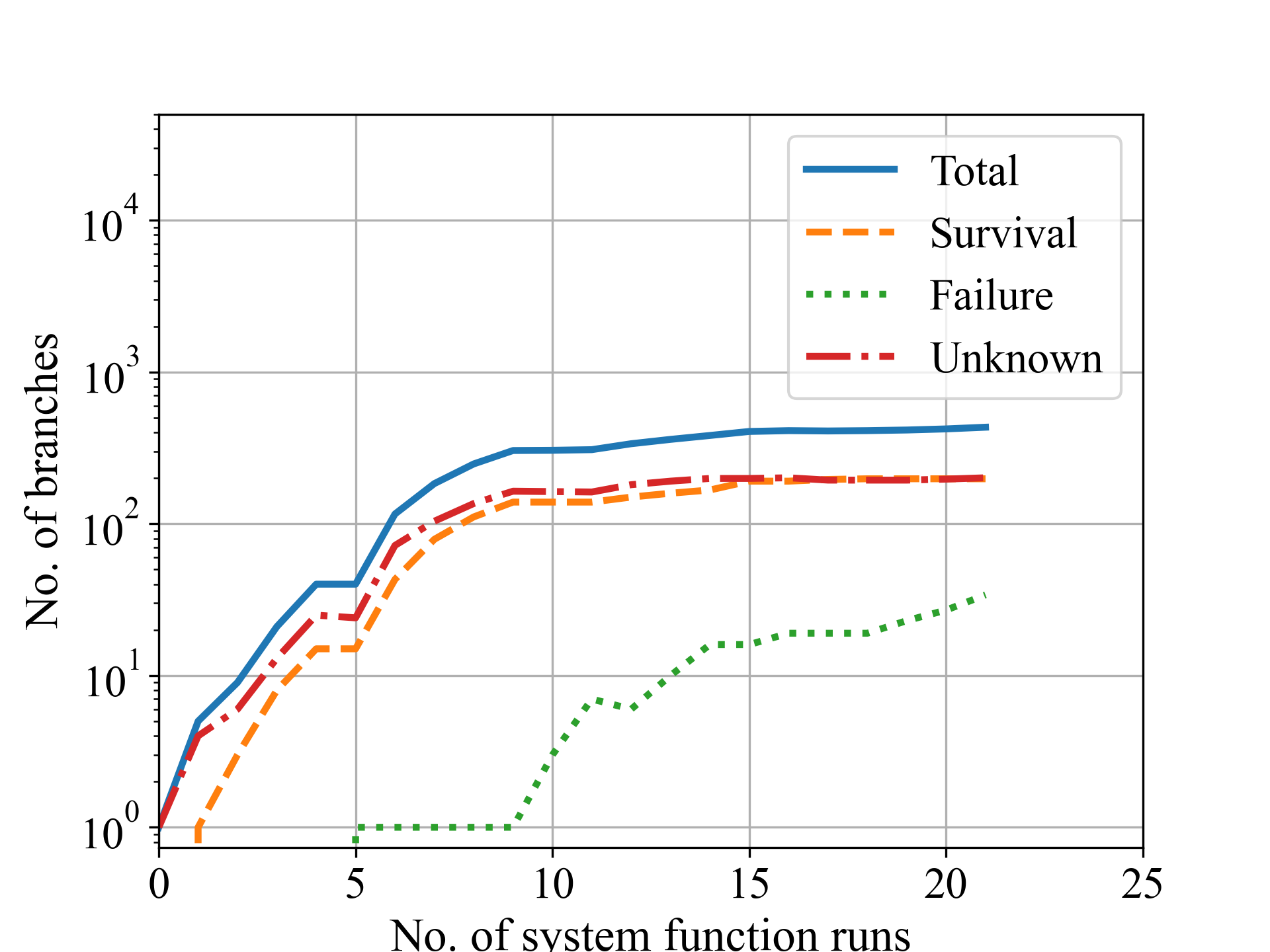}
        \caption{Number of branches (case $d=1,2,3$).}
        \label{fig:br_d1}
    \end{subfigure}
    \hfill
    \begin{subfigure}[b]{0.48\textwidth}
    \centering
        \includegraphics[width=\textwidth]{ 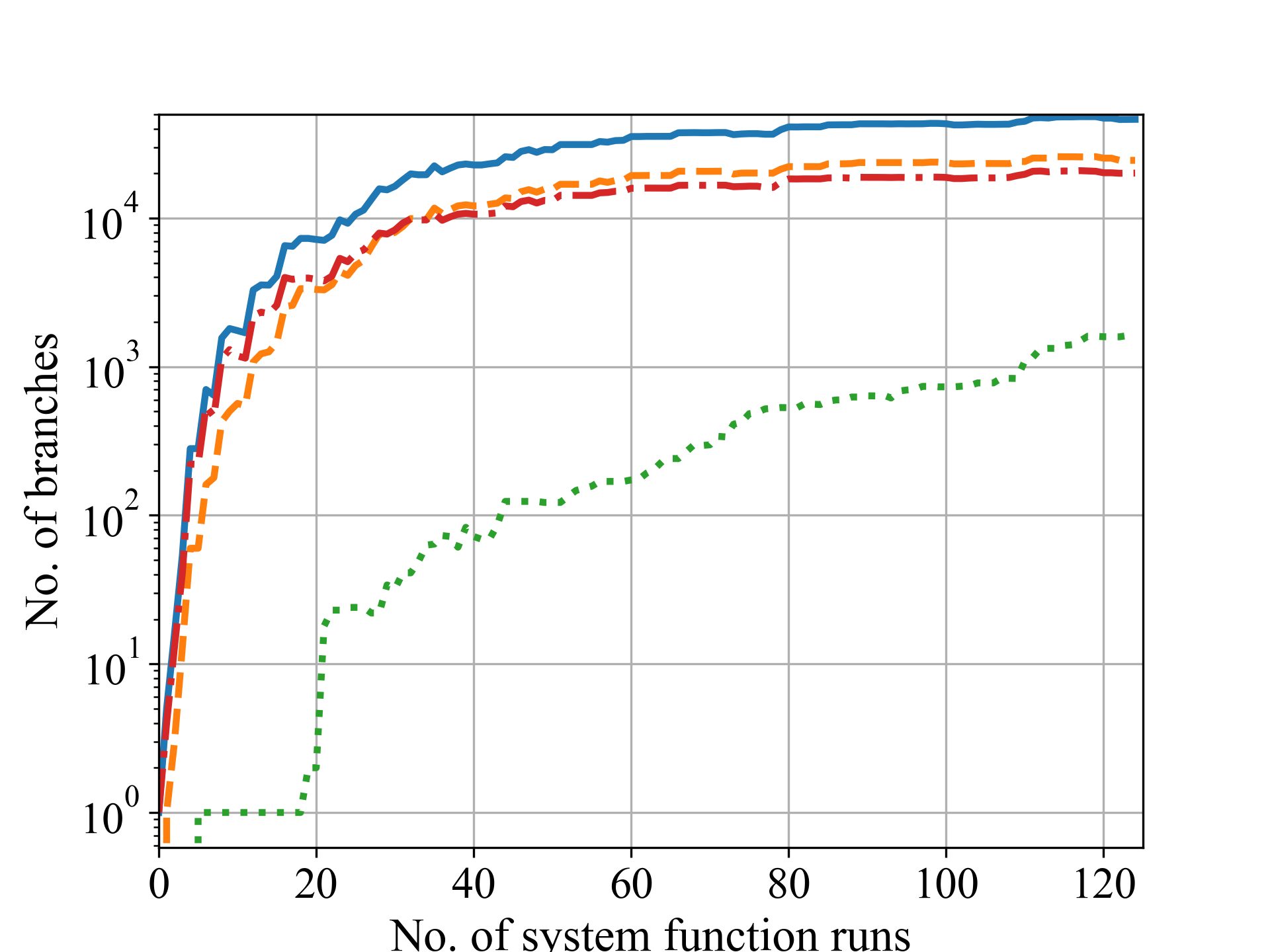}
        \caption{Number of branches (case $d=4,5$).}
        \label{fig:br_d4}
    \end{subfigure}
    \caption{Results of the BRC algorithm for the example network in Fig.~\ref{fig:mf_net}.}
    \label{fig:mf_brc}
\end{figure}

The BRC algorithm is run until a 5\% relative bound width is obtained on the system failure probability. We find that by the construction of the example, the results are the same between the cases $d=1,2,3$ and between the cases $d=4,5$. The convergence of the bounds is illustrated in Fig.~\ref{fig:pf_d1} and \ref{fig:pf_d4} against the number of system function runs, for cases $d=1,2,3$ and $d=4,5$, respectively. The bounds are obtained as $[2.39\cdot 10^{-2}, 2.50\cdot 10^{-2}]$ and $[5.03\cdot 10^{-2}, 5.28\cdot 10^{-2}]$, respectively. This agrees with the Monte Carlo Simulation (MCS) results obtained as $2.49\cdot10^{-2}$ and $5.20\cdot 10^{-2}$ with a coefficient of variance (c.o.v.) 0.01. MCS requires 391,128 and 182,377 samples to reach the c.o.v. of 0.01. The c.o.v. yields the 99\% highest posterior density credible interval of $[2.37, 2.50] \cdot 10^{-2}$ and $[4.82, 5.07] \cdot 10^{-2}$, respectively, where the bound width is 5.3\% of the lower bound \cite{BetPapStr22}. On the other hand, the BRC algorithm requires 22 and 125 runs of the system function to achieve a {\it deterministic} 5\% bound width.

In Fig.~\ref{fig:ru_d1} and \ref{fig:ru_d4}, the number of identified non-dominated rules are plotted. In both cases, the number increases linearly with the number of system function runs, with a final number of 9 (12) and 63 (61) failure (survival) rules. Meanwhile, as illustrated in  Fig.~\ref{fig:br_d1} and Fig.~\ref{fig:br_d4}, the number of branches increases rapidly with the number of system function runs (i.e. the number of rules). We note that the $y$-axes of the figures are in a log-scale. The final numbers of failure, survival, and unspecified branches are 34, 198, 201, and 1,639, 24,481, 20,087, respectively for cases $d=1,2,3$ and $d=4,5$.

In Table~\ref{tab:mf_rules}, the five shortest (i.e. likely to have the highest probabilities) rules are summarised. The short lengths of these rules, from one to five, demonstrate that the BRC algorithm effectively identifies minimal rules. The economical representation of failure rules is notable considering that the system function does not provide any information about failure rules.

In Table~\ref{tab:mf_com}, we compare computational time and the number of system function runs between the BRC algorithm and the complete search by \cite{JanLai08}. The result shows that the BRC algorithm takes much longer computational time and much less system function runs. This is because the primary objective of the BRC algorithm is to minimise the number of system function runs to be able to handle expensive system functions. This is evident by the significantly lower number of system function runs. By construction, the BRC algorithm takes more time than \cite{JanLai08} for decomposition because of the repeated decomposition of every instance a new rule is found. As a result, computational time is significantly longer for the BRC algorithm since a maximum flow analysis requires a very low computational cost taking only fractional seconds. For applications with costly system performance functions, the BRC algorithm would outperform in terms of computational time as well. 

\begin{table}[h!]
    \centering
    \begin{tabular}{c|c|c|c}
    \hline
        \multicolumn{2}{c}{Shortest failure rules} & \multicolumn{2}{|c}{Shortest survival rules} \\ \hline
        $d=1,2,3$ & $d=4,5$ & $d=1,2,3$ & $d=4,5$ \\ \hline
        $(x_{17}^0)$ & $(x_{17}^1)$ & $(x_{3}^1, x_{9}^1, x_{14}^1, x_{17}^1)$ & $(x_{3}^2, x_{9}^2, x_{14}^2, x_{17}^2)$ \\
        $(x_{13}^0, x_{14}^0, x_{16}^0)$ & $(x_{13}^0, x_{14}^0, x_{16}^1)$ & $(x_{2}^1, x_{7}^1, x_{14}^1, x_{17}^1)$ & $(x_{2}^2, x_{7}^2, x_{14}^2, x_{17}^2)$ \\
        $(x_{2}^0, x_{3}^0, x_{4}^0, x_{5}^0)$ & $(x_{13}^0, x_{14}^0, x_{18}^1)$ & $(x_{4}^1, x_{13}^1, x_{17}^1, x_{20}^1)$ & $(x_{4}^2, x_{13}^2, x_{17}^2, x_{20}^2)$ \\
        $(x_{5}^0, x_{13}^0, x_{14}^0, x_{21}^0)$ & $(x_{13}^0, x_{14}^1, x_{16}^0)$ & $(x_{2}^1, x_{16}^1, x_{17}^1, x_{18}^1, x_{21}^1)$ & $(x_{1}^2, x_{5}^2, x_{16}^2, x_{17}^2, x_{18}^2)$ \\
        $(x_{1}^0, x_{2}^0, x_{3}^0, x_{4}^0)$ & $(x_{4}^0, x_{14}^0, x_{16}^1)$ & $(x_{1}^1, x_{5}^1, x_{16}^1, x_{17}^1, x_{18}^1)$ & $(x_{4}^2, x_{9}^2, x_{11}^2, x_{14}^2, x_{17}^2)$ \\ \hline
        
    \end{tabular}
    \caption{The five shortest failure and survival rules identified by the BRC algorithm.}
    \label{tab:mf_rules}
\end{table}

\begin{table}[h!]
    \centering
    \begin{tabular}{c|c|c|c|c}
    \hline
        & \multicolumn{2}{c|}{BRC algorithm} & \multicolumn{2}{c}{Complete search by  \cite{JanLai08}} \\ \hline
        Failure criterion & $d=1,2,3$ & $d=4,5$ & $d=1,2,3$ & $d=4,5$ \\ \hline
        \makecell{No. of system runs} & 22 & 125 & 1,272$^\text{a}$ & 40,004 \\
        \makecell{Computation time \\[-5pt] (sec.)} & 2.28 & 2,831 & 0.016 & 0.268 \\
        \hline
        \multicolumn{5}{l}{\small $^\text{a}$ In \cite{JanLai08}, a system function is run twice during a while loop, each on a lower and an upper bound.}
        
    \end{tabular}
    \caption{Comparison of computational performance between the BRC algorithm and a complete search by \cite{JanLai08}.}
    \label{tab:mf_com}
\end{table}

\subsection{Eastern Massachusetts highway benchmark network: reliability mapping and updating}\label{sec:ema}

\subsubsection{Settings}\label{sec:ema_set}
\begin{figure}[h!]
    \centering
    \includegraphics[width=0.7\textwidth]{ 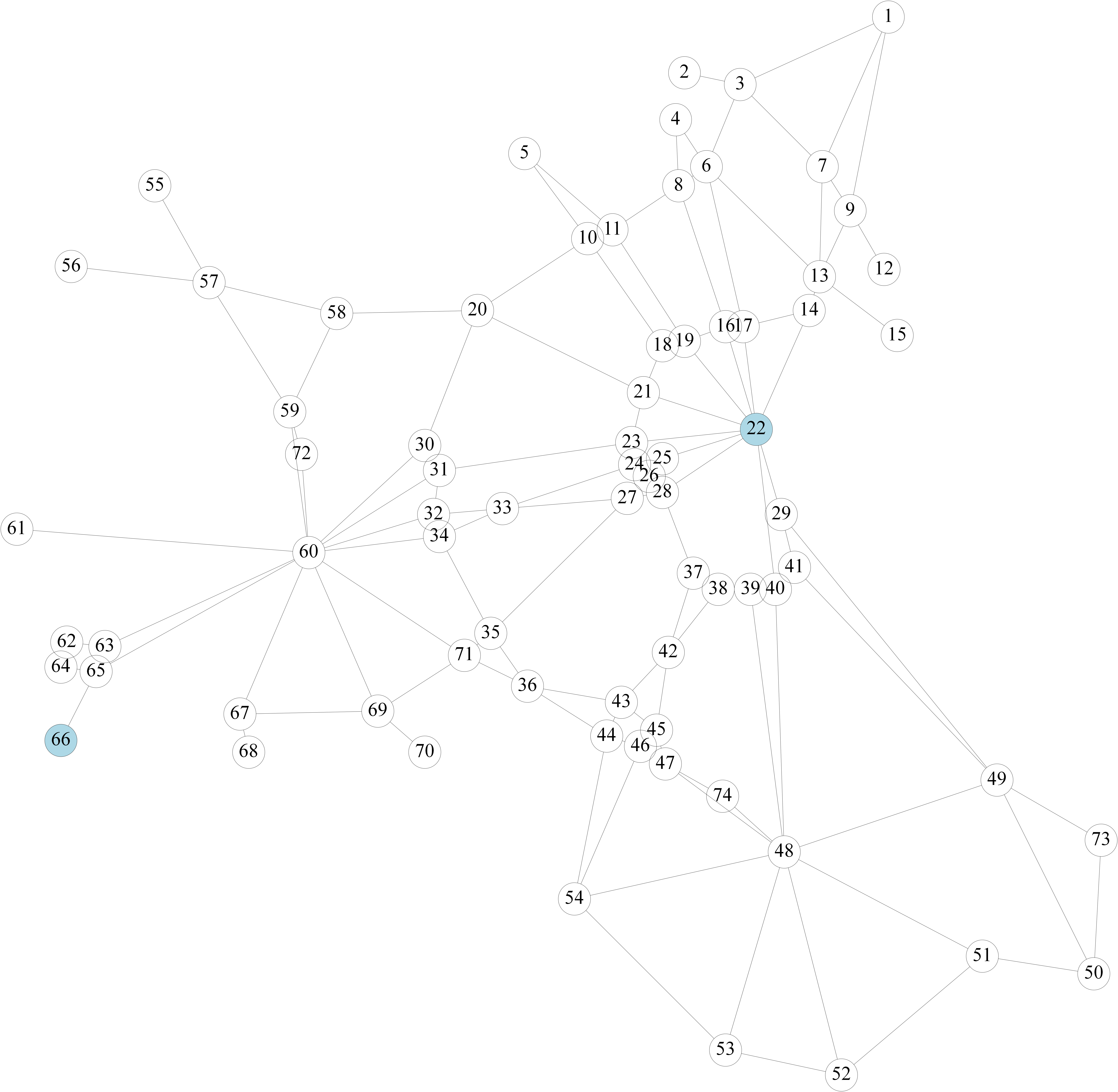}
    \caption{EMA highway benchmark network analysed in Section~\ref{sec:ema}. Origin nodes are marked by light blue. Each of the non-origin nodes constitutes a system event, whose failure event is defined as the distance to the closest origin being longer than twice a pre-disaster distance.}
    \label{fig:ema}
\end{figure}

To study the performance of the BRC algorithm on a more realistic and larger system, we investigate the Eastern Massachusetts (EMA) highway benchmark network presented in Fig.~\ref{fig:ema}. The network consists of 74 nodes and 129 edges, which represent roadways. The connectivity information is summarised in Table~\ref{tab:ema_edge}. We consider seismic risks as a potential cause for the failure of roadways. The state of each edge $e_1,\ldots,e_{129}$ is represented by a random variable $X_1,\ldots,X_{129}$, which takes a binary state, i.e. 0 for failure and 1 for survival. The probability distribution $P(\bm{X})$ is specified in Section~\ref{sec:ema_fra}.

To define system failure, we set international airports as the origin nodes since after an earthquake, they are used as distribution hubs of emergency resources. Two airports are identified at nodes $n_{22}$ and $n_{66}$, which are marked by light blue in Fig.~\ref{fig:ema}. For any node, system failure is defined that the node's distance to the closest origin node is longer than twice the pre-disaster shortest distance. Since there are 72 non-origin nodes, this definition leads to 72 system failure events. The following analyses aim to generate a map of failure probabilities over the network and update the failure probabilities given new hazard information, which requires evaluating the probability of these 72 system failure events.

The system performance function evaluates the shortest travel path and time to each of the two origin nodes. Then, it selects the shorter path to determine the system state. In case of system survival, the system function returns a minimal survival rule as the states of the component events constituting a selected shortest path. For instance, for $n_{31}$, if the shortest path algorithm returns a shortest path to $n_{66}$ as $e_{62}\rightarrow e_{117} \rightarrow e_{125}$ (corresponding to a node path $n_{31}\rightarrow n_{60} \rightarrow n_{65} \rightarrow n_{66}$), one can identify $(x_{62}^1,x_{117}^1,x_{125}^1)$ as a corresponding survival rule. By contrast, in case of system failure, the system function does not return any minimal failure rule. For computation, we use a CPU Intel{\textregistered} Core{\texttrademark} i9 and a 64 GB RAM.

\subsubsection{Probability of component system states}\label{sec:ema_fra}

To model earthquake loading on the components, we use the ground motion propagation equation (GMPE) developed for Eastern North America (ENA) by \cite{Cam03}, given as
\begin{equation}\label{eq:gmpe}
    \ln Y = c_1 + f_1(M_W) + f_2 (M_W, r_{\text{rup}}) + f_3(r_{\text{rup}}),
\end{equation}
where 
\begin{subequations}\label{eq:gmpe_part2}
    \begin{equation}
        f_1(M_W) = c_2 M_W + c_3(8.5-M_W)^2,
    \end{equation}
    \begin{equation}
        f_2(M_W, r_{\text{rup}}) = c_4 \ln R + (c_5 + c_6 M_W) r_{\text{rup}},
    \end{equation}
    \begin{equation}
        R = \sqrt{r_{\text{rup}}^2 + [c_7 \exp (c_8 M_W)]^2}, \  \text{and}
    \end{equation}
    \begin{equation}
        f_3(r_{\text{rup}}) = \begin{cases}
            0 & \text{for $r_{\text{rup}} \leq r_1$} \\
            c_7(\ln r_{\text{rup}} - \ln r_1) & \text{for $r_1 < r_{\text{rup}} \leq r_2$} \\
            c_7(\ln r_{\text{rup}} - \ln r_1) + c_8(\ln r_{\text{rup}} - \ln r_2) & \text{for $r_{\text{rup}} > r_2$}.
        \end{cases}
    \end{equation}
\end{subequations}
We assume that all highway segments include a bridge, and bridge failure is the primary cause of road closure. Therefore, in equations \eqref{eq:gmpe} and \eqref{eq:gmpe_part2}, $Y$ is 5\% damped pseudo spectral acceleration (PSA) in gravitational acceleration (g), $M_W$ is moment magnitude, and $r_{\text{rup}}$ is the closest distance to fault rupture (km), and $r_1$ and $r_2$ are given as 70 km and 130 km, respectively. The epistemic standard deviation of $\ln Y$ is quantified as 
\begin{equation}\label{eq:AleStaDev}
    \sigma_{\ln Y} = \begin{cases}
        c_{11} + c_{12}M_W & \text{for $M_W < M_1$} \\
        c_{13} & \text{for $M_W \geq M_1$}.
    \end{cases}
\end{equation}
The coefficients for natural period 1.0 second are $c_1=-0.6104$, $c_2 = 0.451$, $c_3 = -0.2090$, $c_4 = -1.158$, $c_5 = -0.00255$, $c_6 = 0.000141$, $c_7 = 0.299$, $c_8=0.503$, $c_{11} = 1.110$, $c_{12}=-0.0793$, $c_{13}=0.543$, and $M_1 = 7.16$. The epistemic standard deviation $\sigma_{\ln Y}^e$ is provided in \ref{tab:EpiStd}, which is excerpted from \citep{Cam03}.

For fragility modelling, we hypothetically classify roadways into three fragility classifications by HAZUS as summarised in Table~\ref{tab:ema_fra}. We assume that a component fails, i.e. $X_n = 0$, if an edge $e_n$ has a damage state equal to or worse than an extensive damage state. To calculate probabilities $P(x_n^0)$, we assume that seismic demand on $e_n$, i.e. $Y_n$ in \eqref{eq:gmpe}, and its structural capacity, i.e. $R_n$ from Table~\ref{tab:ema_fra}, are statistically independent. Also, we consider epistemic uncertainty $\sigma_{\ln Y}^e$ (and not aleatory uncertainty $\sigma_{\ln Y}$) as a source of correlations in ground motions. In other words, a joint probability is calculated as
\begin{equation}
   P(x_{n_1}, \ldots, x_{n_M}) = \Phi( \textbf{0}; \bm{\mu}, \mathbf{\Sigma} ),
\end{equation}
where \textbf{0} is a zero-vector with length $M$, and $\Phi(\cdot)$ is the cumulative distribution function of a multivariate normal distribution with mean $\bm{\mu}=(\mu_{n_1}, \ldots, \mu_{n_M})$ and covariance $\mathbf{\Sigma}$ such that
\begin{subequations}
    \begin{equation}
        \mu_n = \begin{cases}
            -\ln R_n + \ln Y_n & \text{for $X_{n}=1$} \\
            \ln R_n - \ln Y_n & \text{for $X_n = 0$},
        \end{cases}
    \end{equation}
    \begin{equation}\label{eq:cov_ema}
        \Sigma_{nm} = \begin{cases}
            \sigma_{\ln Y}^2 + (\sigma_{\ln Y_n}^{e})^2 + \beta_n^2 & \text{for $n=m$} \\
            \sigma_{\ln Y_n}^{e} \cdot \sigma_{\ln Y_m}^{e} & \text{for $n\neq m$}.
        \end{cases}
    \end{equation}
\end{subequations}
In the equations above, $R_n$ and $\beta_n$ denote respectively the logarithm of the mean fragility value and the dispersion summarised in Table~\ref{tab:ema_fra}; $\ln Y_n$ is the mean seismic demand on $e_n$ calculated by \eqref{eq:gmpe}; and $\sigma_{\ln Y}$ and $\sigma_{\ln Y_n}^{e}$ represent respectively the aleatory and the epistemic uncertainties in the seismic demand (i.e. \eqref{eq:AleStaDev} and Table~\ref{tab:EpiStd}).

\subsubsection{Results}\label{sec:ema_res}

\begin{figure}[h!]
    \centering
    \begin{subfigure}[b]{0.48\textwidth}
        \centering
        \includegraphics[height=0.75\textwidth]{ 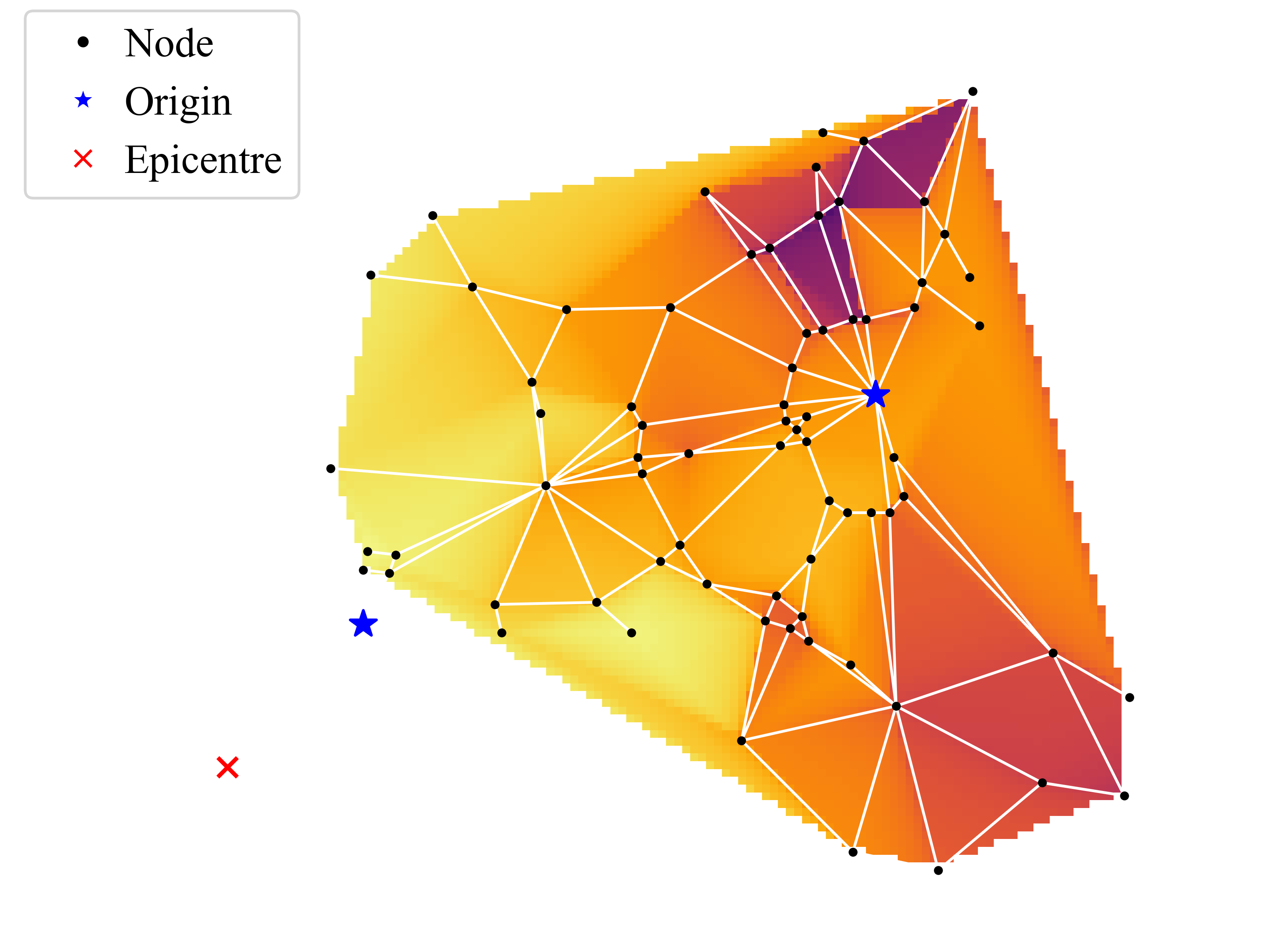}
        \caption{Mw = 8.0. Component events are considered statistically independent.}
        \label{fig:map_ind}
    \end{subfigure}    
    \hfill
    \begin{subfigure}[b]{0.48\textwidth}
        \centering
        \includegraphics[height=0.75\textwidth]{ 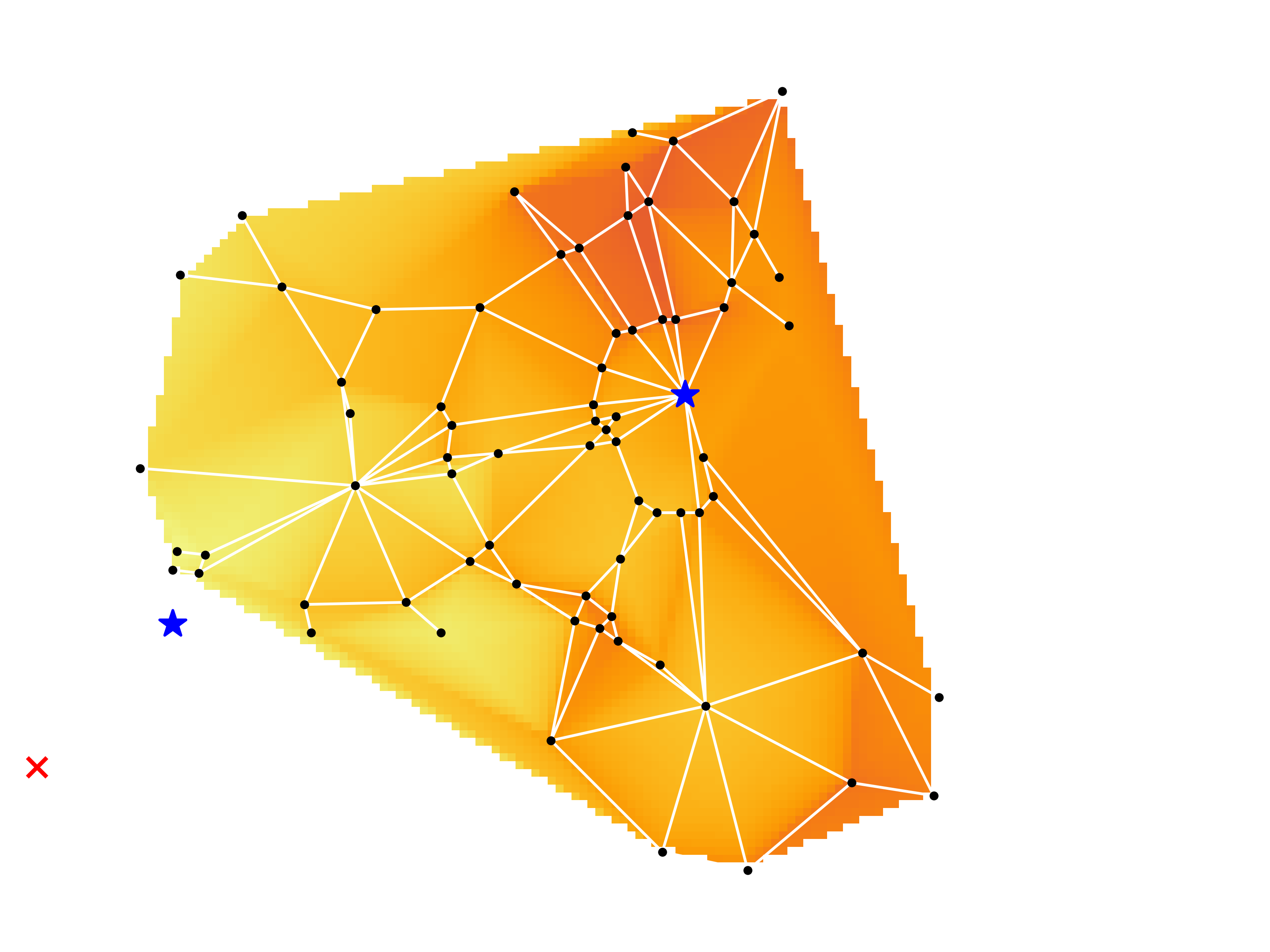}
        \caption{Mw = 8.0. Component events are considered statistically dependent.}
        \label{fig:map_dep}
    \end{subfigure}
    
    \begin{subfigure}[b]{0.48\textwidth}
        \centering
        \includegraphics[height=0.75\textwidth]{ 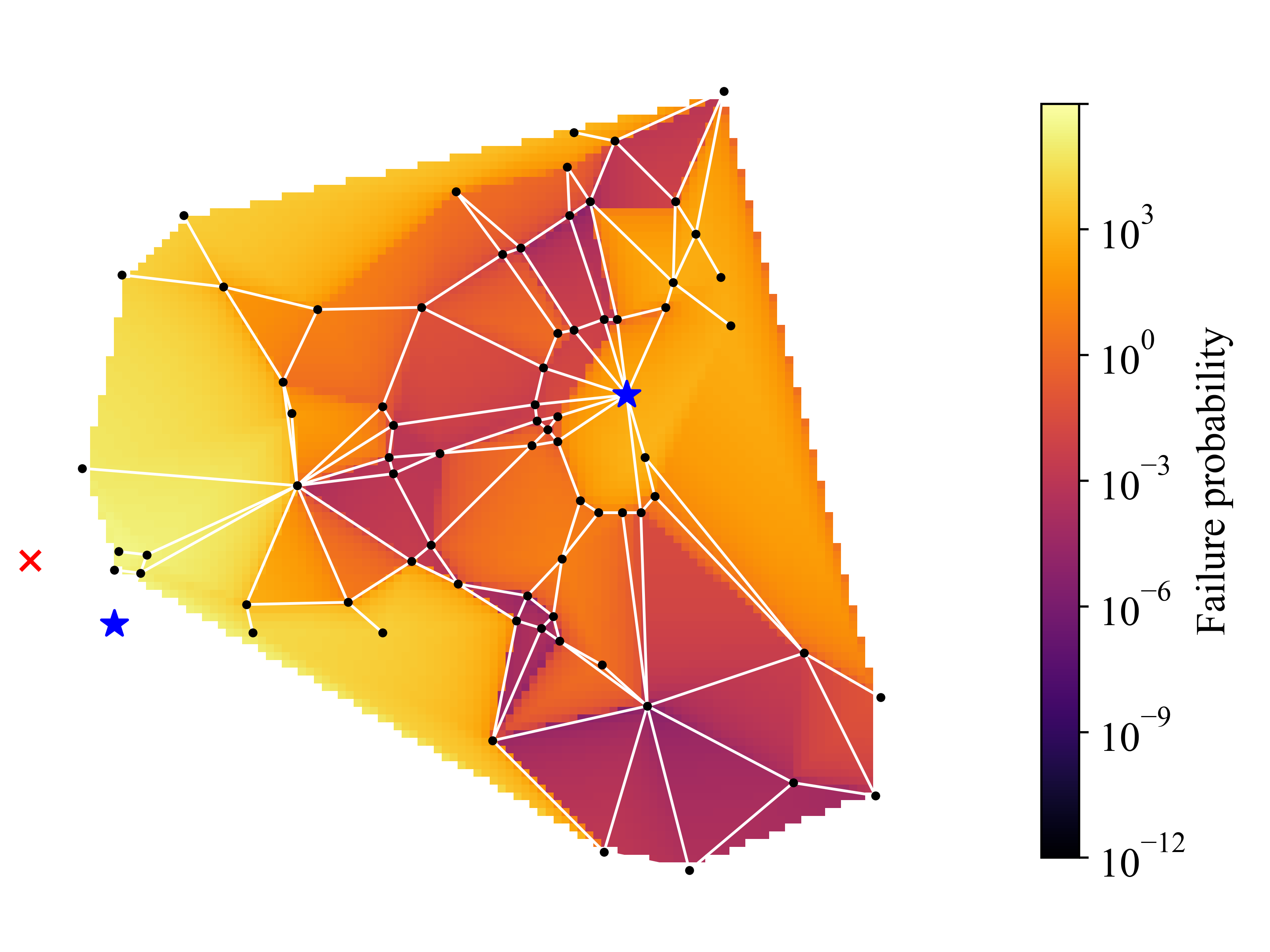}
        \caption{Updated by a new hazard scenario with Mw = 8.2. Component events are considered statistically independent.}
        \label{fig:map_upd}
    \end{subfigure}

    \caption{Linearly interpolated maps of the logarithm of failure probabilities of the non-origin nodes in the EMA benchmark network. The two origin nodes and epicentre locations are marked in the figure by blue stars and red crosses, respectively.}
    \label{fig:map}
\end{figure}

We run the BRC algorithm by assuming an epicenter with magnitude $M_w=8.0$ and location as marked in Fig.~\ref{fig:map_ind}. The BRC algorithm is performed until either 5\% of a bound width is reached, or the total number of (failure, survival, and unspecified) branches reaches 50,000. In the latter case, an additional MCS is performed on unspecified branches until the c.o.v. is less than 0.01.
The obtained failure probabilities (more precisely, the lower bound of either the 5\%-bound on the failure probability or the 95\%-credible interval) for each node are visualised in Fig.~\ref{fig:map_ind}, where spaces between nodes are interpolated linearly with respect to the logarithm of the failure probabilities. After applying the BRC algorithm, probabilities of the identified branches are re-computed considering the dependence between component events, which is presented in Fig.~\ref{fig:map_dep}. In addition, assuming that an earthquake has occurred and the map needs to be updated in near-real-time to inform emergency activities, we update the probabilities of the branches and the samples assuming a hazard scenario with magnitude $M_w = 8.2$ and epicentre location as illustrated in Fig.~\ref{fig:map_upd}. In this case, we ignore dependence between components by assuming that information on actual ground motions is available and accurate enough. The three results show how the results change by the consideration of dependence and the new hazard scenario. For validation, we compare the results of selected nodes with MCS results in Table~\ref{tab:pfs}; as can be seen, the two results agree.

\begin{figure}[h!]
    \centering
    \begin{subfigure}[b]{\textwidth}
        \includegraphics[width=\textwidth]{ 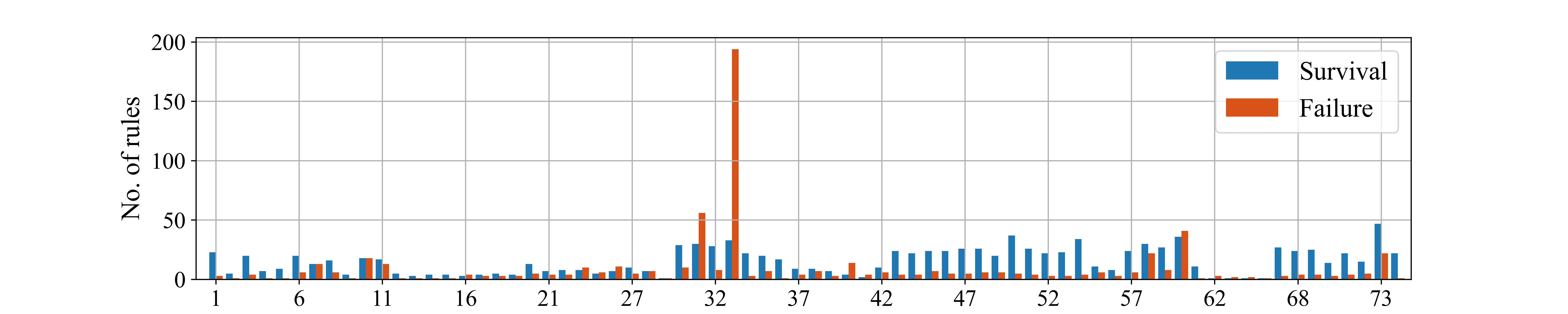}
        \caption{The number of identified non-dominated failure and survival rules.}
        \label{fig:n_rule}
    \end{subfigure}
    \begin{subfigure}[b]{0.89\textwidth}
        \centering
    \includegraphics[width=\textwidth]{ 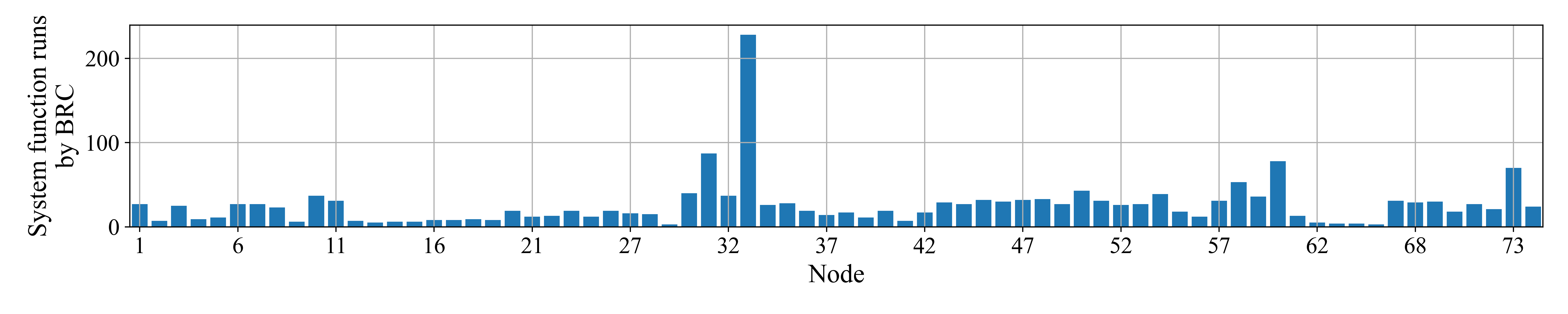}
    \caption{The number of system function runs to obtain the non-dominated rules.}
    \label{fig:n_funs_brc}
    \end{subfigure}
    
    \begin{subfigure}[b]{0.6\textwidth}
        \centering
        \includegraphics[width=0.5\textwidth]{ 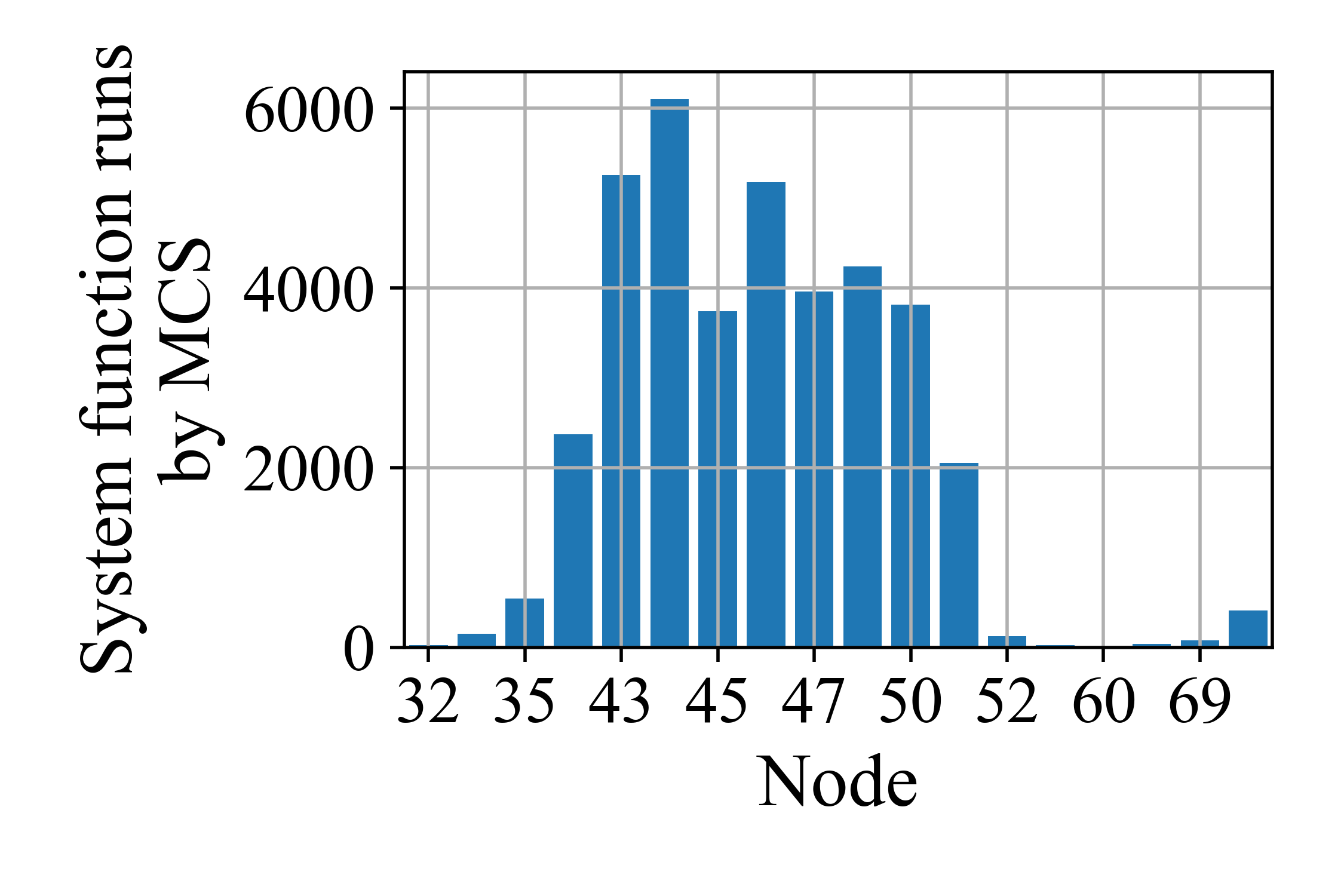}
        \caption{The number of MCS samples on unspecified branches. MCS is performed until a c.o.v. is less than 0.01. MCS is performed on nodes $n_{32}$, $n_{34}$-$n_{36}$, $n_{43}$, $n_{44}$-$n_{48}$, $n_{50}$-$n_{52}$, $n_{54}$, $n_{60}$, $n_{67}$, $n_{69}$, $n_{71}$.}
        \label{fig:n_funs_mcs}
    \end{subfigure}
\caption{The number of system function runs to obtain the results in Fig.~\ref{fig:map}. The $x$-axes represent the non-origin nodes.}
\end{figure}

The BRC algorithm achieves a 5\%-bound on the system failure probability for 54 nodes by decomposition only; and it performs hybrid inference for the remaining 18 nodes. The number of identified non-dominated rules are presented in Fig.~\ref{fig:n_rule}, where all nodes have less than 100 rules in total except for $n_{33}$ which has 227 rules. This shows that despite the large size of the network, systems can be represented by a tractable number of rules. In Fig.~\ref{fig:n_funs_brc}, the number of system function runs required by the BRC algorithm is summarised, where again all numbers are less than 100 except for $n_{33}$ with 228. In other words, the BRC algorithm finds minimal rules effectively almost every time a system performance function is run. As summarised in Fig.~\ref{fig:n_funs_mcs}, for the 18 nodes MCS is applied, the maximum number of MCS samples to achieve c.o.v. 0.01\% is 6,102 for $n_{44}$. Overall, all system events take less than 6,129 number of system function runs, which is very low considering that some of the events have very low failure probabilities on the order of $10^{-11}$ (for $n_6$). 

For illustration, Table~\ref{tab:ema_rules} presents the five shortest failure and survival rules for $n_{30}$ and $n_{73}$; a total of 22 (47) and 10 (29) failure (survival) rules have been identified for these two nodes, respectively. The listed survival rules show the shortest, and likely most reliable, paths to the origin nodes. This underlines an additional benefit of the BRC algorithm: it identifies critical link- and cut-sets of a system event, which provides useful information for risk management. For example, in Table~\ref{tab:ema_rules}, the first two survival rules for $n_{30}$ connect the node to the origin $n_{22}$ and the last three to $n_{66}$. During an emergency situation, provided real-time updates on damaged roadways, such information can be used to decide on an emergency route to the area.

\begin{table}[h!]
    \centering
    \begin{tabular}{c|c|c|c}
    \hline
        \multicolumn{2}{c}{Shortest failure rules} & \multicolumn{2}{|c}{Shortest survival rules} \\ \hline
        $n_{30}$ & $n_{73}$ & $n_{30}$ & $n_{73}$ \\ \hline
        $(x_{36}^0, x_{59}^0, x_{60}^0)$ & $(x_{102}^0, x_{104}^0)$ & $(x_{40}^1, x_{46}^1, x_{59}^1)$ & $(x_{43}^1, x_{58}^1, x_{102}^1)$ \\
        $(x_{22}^0, x_{35}^0, x_{59}^0, x_{60}^0)$ & $(x_{101}^0, x_{102}^0, x_{103}^0)$ & $(x_{35}^1, x_{36}^1, x_{38}^1)$ & $(x_{43}^1, x_{57}^1, x_{82}^1, x_{102}^1)$ \\
        $(x_{36}^0, x_{46}^0, x_{60}^0, x_{61}^0, x_{62}^0)$ & $(x_{43}^0, x_{44}^0, x_{55}^0, x_{56}^0)$ & $(x_{60}^1, x_{117}^1, x_{125}^1)$ & $(x_{44}^1, x_{80}^1, x_{82}^1, x_{102}^1)$ \\
        $(x_{36}^0, x_{59}^0, x_{62}^0, x_{65}^0, x_{68}^0, x_{125}^0)$ & $(x_{43}^0, x_{44}^0, x_{56}^0, x_{95}^0)$ & $(x_{60}^1, x_{116}^1, x_{123}^1, x_{125}^1)$ & $(x_{43}^1, x_{58}^1, x_{101}^1, x_{104}^1)$ \\
        $(x_{22}^0, x_{35}^0, x_{46}^0, x_{60}^0, x_{61}^0, x_{62}^0)$ & $(x_{96}^0, x_{97}^0, x_{101}^0, x_{102}^0)$& $(x_{59}^1, x_{62}^1, x_{117}^1, x_{125}^1)$ & $(x_{44}^1, x_{81}^1, x_{95}^1, x_{102}^1)$ \\ \hline
        
    \end{tabular}
    \caption{The five shortest failure and survival rules of $n_{30}$ and $n_{73}$. For $n_{30}$ and $n_{73}$, 22 (47) and 10 (29) of non-dominated failure (survival) rules are identified, respectively. With respect to independent component distributions, bounds on failure probability are computed as $[1.32, 1.38]\cdot 10^{-3}$ and $[5.34, 5.61]\cdot 10^{-6}$ (for both system events, sampling is not performed because of the sufficiently narrow bound width).}
    \label{tab:ema_rules}
\end{table}

To assess computational efficiency, we arbitrarily select six nodes in Table~\ref{tab:pfs} and summarise their failure probabilities and computational time taken by the BRC algorithm and MCS. The nodes have a varying order of failure probability from $3.33\cdot 10^{-3}$ to $3.92\cdot 10^{-1}$. Nodes with lower probabilities are not selected to avoid excessive computational cost required by MCS (obviously, BRC would outperform MCS even more for those nodes). We confirm that in all cases the BRC algorithm takes a shorter time than MCS: it takes from $3\cdot 10^{-5}$ \% (for node $n_{29}$) to $79$ \% ($n_{67}$) of the time taken by MCS. We also compare the computational time in the case of probability updating, in which case the reduction is much more significant: it takes at maximum $0.3$ \% (for $n_{64}$) of the time taken by MCS.

\begin{figure}[h!]
    \centering
    \begin{subfigure}[b]{\textwidth}
        \centering
        \includegraphics[width=0.9\textwidth]{ 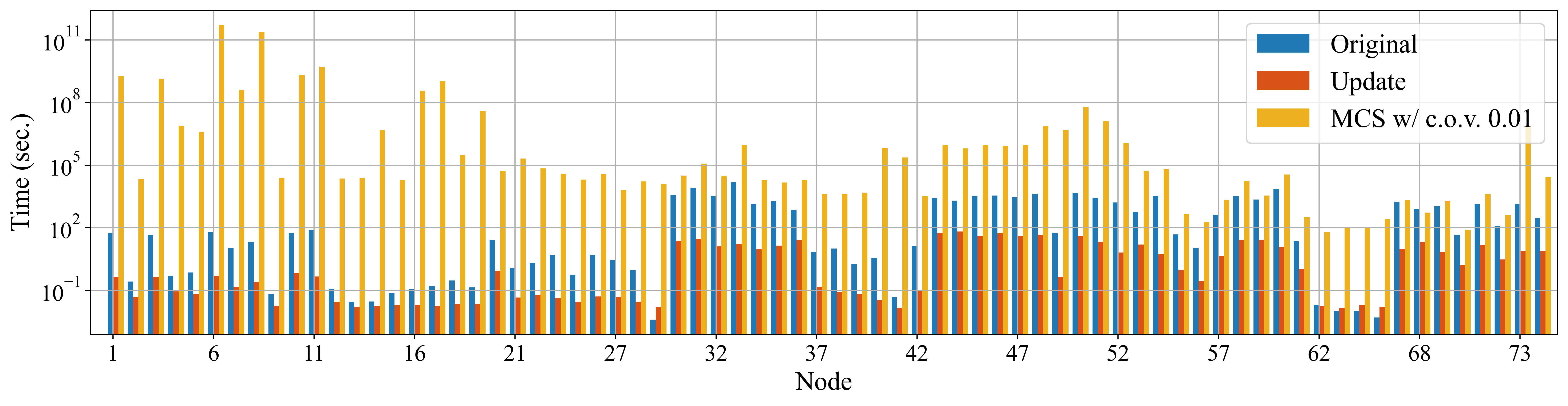}
        \caption{Computational time (sec.) taken for analysis. Referring to the results in Table~\ref{tab:pfs}, the time taken for MCS is calculated by assuming 0.00395 seconds for each run of a system performance function. The total time taken for all of 72 nodes are 24.3 hours, 9.88 minutes, and $2.13\cdot 10^8$ hours respectively for the original analysis (including the BRC algorithm and MCS on unspecified branches for c.o.v. 0.01), for probability update, and MCS for c.o.v. 0.01.}
        \label{fig:time}
    \end{subfigure}

    \begin{subfigure}[b]{\textwidth}
        \centering
        \includegraphics[width=\textwidth]{ 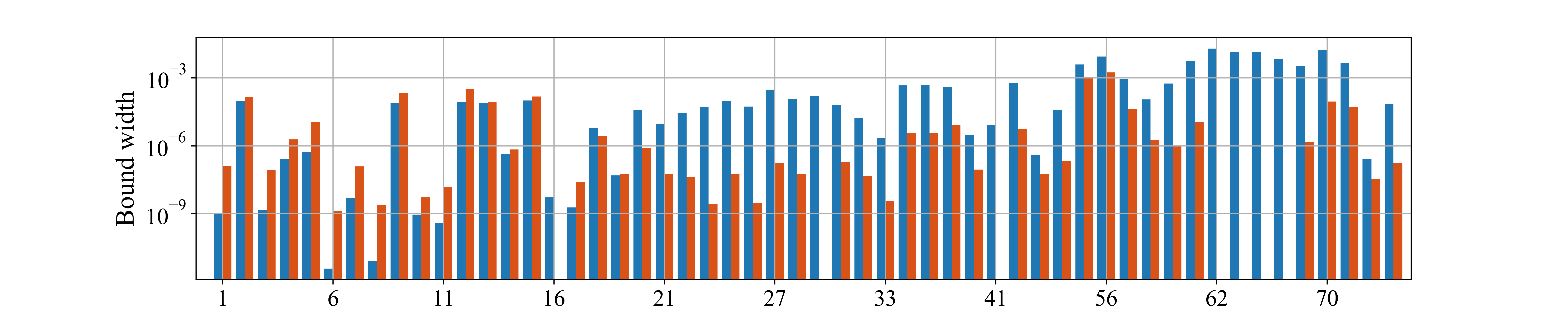}
        \caption{Bound widths for the nodes where only the BRC algorithm is applied.}
        \label{fig:bnd_wid}
    \end{subfigure}

    \begin{subfigure}[b]{0.9\textwidth}
        \centering
        \includegraphics[scale=0.5]{ 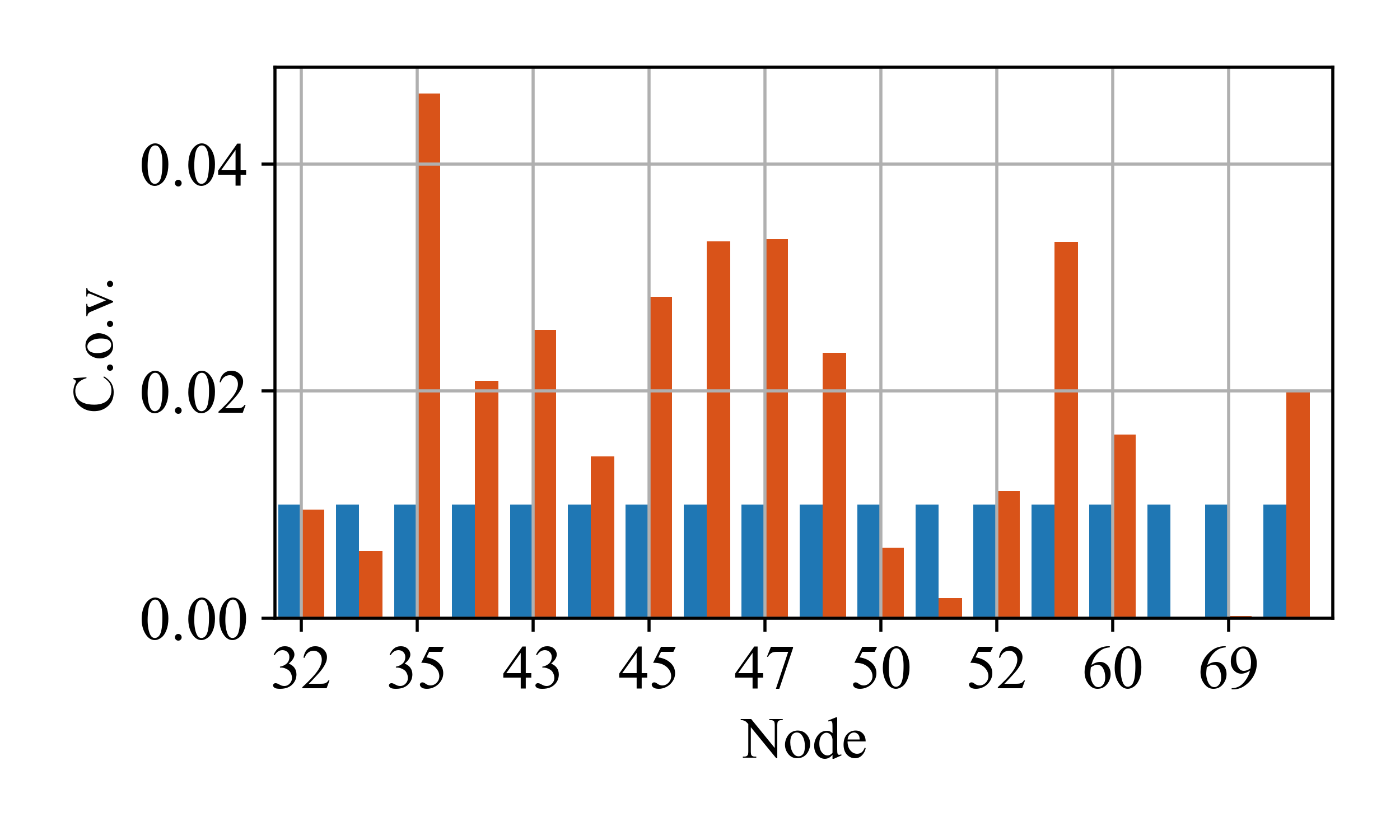}
        \caption{C.o.v.'s for the nodes where hybrid inference is applied.}
        \label{fig:cov}
    \end{subfigure}
\caption{Comparison between the original decomposition by the BRC algorithm, the updated result by a new hazard scenario, and MCS.}
\end{figure}

We compare computational times in Fig.~\ref{fig:time} for the initial run of the BRC algorithm (to obtain Fig.~\ref{fig:map_ind}), the probability update (Fig.~\ref{fig:map_upd}), and, for comparison, MCS for c.o.v. 0.01 given the initial setting of Fig.~\ref{fig:map_ind}. Since MCS takes an unaffordable computational time for many nodes, we evaluate the computational times by assuming 0.00395 seconds for a system function run, which is the average time taken for the analyses of Table~\ref{tab:pfs}; the required number of samples is evaluated based on the evaluated failure probabilities by the BRC algorithm. As illustrated in Fig.~\ref{fig:time}, the computation times of the BRC algorithm (blue and red bars) are orders of magnitudes lower than those of MCS (yellow bars) for most cases. The total time taken for all of the 72 system events is 24.3 hours, the time taken by BRC for updating probabilities is 9.88 minutes, and the time taken for MCS is $2.13\cdot 10^8$ hours. The particularly short time for the probability update indicates that for risk management, one can apply the BRC algorithm beforehand and facing a hazard event quickly update the probabilities to coordinate emergency activities. We note that the computation can be multi-processed for each node to further reduce the computation time. 

To check if a probability update does not compromise analysis precision, we compare bound widths of failure probability for nodes with only the BRC algorithm applied on in Fig.~\ref{fig:bnd_wid} and c.o.v.'s for those with hybrid inference applied on in Fig.~\ref{fig:cov}, respectively. We note that all bounds remain within the same order, while the maximum c.o.v. (observed for $n_{35}$) is less than 0.05. We even observe that in some cases, bound width as well as c.o.v. are even less than original results. For instance, in Fig.~\ref{fig:cov}, $n_{67}$ has c.o.v. $9.96\cdot 10^{-3}$ originally, but with the updated probabilities, has $3.06\cdot 10^{-6}$. This is because in the original analysis, the probability of unspecified branches is $2.33\cdot10^{-3}$ with $1.79\cdot 10^{-2}$ of failure branches, while the updated probabilities lead to $1.36\cdot 10^{-6}$ with $6.05\cdot 10^{-3}$ of failure branches. Such well-controlled precision of inference results demonstrates that the the BRC algorithm, as a non-simulation-based method, allows analysis results to remain robust against changing probabilities.

\begin{table}[h!]
    \centering
    \begin{tabular}{c|c c | c c}
        \hline 
        \multirow{2}{*}{Node} & \multicolumn{2}{c|}{Original} & \multicolumn{2}{c}{Updated} \\
        \cline{2-5}
         & \makecell[c]{The BRC algorithm \\[-5pt] (time)} & \makecell[c]{MCS (no. of samples, \\[-5pt] time)} & \makecell[c]{The BRC algorithm \\[-5pt] (time)} & \makecell[c]{MCS (no. of samples, \\[-5pt] time)} \\
         \hline
        $n_{29}$ & \makecell[c]{$3.33\cdot 10^{-3}$ \\[-5pt] ($4.00\cdot 10^{-3}$)} & \makecell[c]{$3.28\cdot 10^{-3}$ ($3.03\cdot 10^6$, \\[-5pt] $1.19\cdot 10^4$)} & \makecell[c]{$9.44\cdot 10^{-3}$ \\[-5pt] ($1.60\cdot 10^{-2}$)} & \makecell[c]{$9.52\cdot 10^{-3}$ ($1.04\cdot 10^6$, \\[-5pt] $4.08\cdot 10^3$)} \\
        $n_{62}$ & \makecell[c]{$3.92\cdot 10^{-1}$ \\[-5pt] ($2.00\cdot 10^{-2}$)} & \makecell[c]{$3.94\cdot 10^{-1}$ ($1.54\cdot 10^4$, \\[-5pt] $6.11\cdot 10^1$)} & \makecell[c]{$4.39\cdot 10^{-1}$ \\[-5pt] ($1.70\cdot 10^{-2}$)} & \makecell[c]{$4.38\cdot 10^{-1}$ ($1.29\cdot 10^4$, \\[-5pt] $4.98\cdot 10^1$)} \\
        $n_{63}$ & \makecell[c]{$2.68\cdot 10^{-1}$ \\[-5pt] ($1.00\cdot 10^{-2}$)} & \makecell[c]{$2.70\cdot 10^{-1}$ ($2.71\cdot 10^4$, \\[-5pt] $1.05\cdot 10^{2}$)} & \makecell[c]{$3.01\cdot 10^{-1}$ \\[-5pt] ($1.40\cdot 10^{-2}$)} & \makecell[c]{$3.02\cdot 10^{-1}$ ($2.31\cdot 10^4$, \\[-5pt] 9.00)} \\
        $n_{64}$ & \makecell[c]{$2.84\cdot 10^{-1}$ \\[-5pt] ($1.00\cdot 10^{-2}$)} & \makecell[c]{$2.84\cdot 10^{-1}$ ($2.53\cdot 10^4$, \\[-5pt] $9.40\cdot 10^1$)} & \makecell[c]{$3.43\cdot 10^{-1}$ \\[-5pt] ($1.90\cdot 10^{-2}$)} & \makecell[c]{$3.42\cdot 10^{-1}$ ($1.93\cdot 10^4$, \\[-5pt] 7.20)} \\
        $n_{65}$ & \makecell[c]{$1.33\cdot 10^{-1}$ \\[-5pt] ($5.00\cdot 10^{-3}$)} & \makecell[c]{$1.32\cdot 10^{-1}$ ($6.56\cdot 10^4$, \\[-5pt] $2.50\cdot 10^{2}$)} & \makecell[c]{$1.77\cdot 10^{-1}$ \\[-5pt] ($1.60\cdot 10^{-2}$)} & \makecell[c]{$1.76\cdot 10^{-1}$ ($4.58\cdot 10^4$, \\[-5pt] $1.78\cdot 10^{2}$)} \\
        $n_{67}$ & \makecell[c]{$1.86\cdot 10^{-2}$ \\[-5pt] ($1.79\cdot 10^{3}$)} & \makecell[c]{$1.81\cdot 10^{-2}$ ($5.42\cdot 10^5$, \\[-5pt] $2.26\cdot 10^{3}$)} & \makecell[c]{$6.05\cdot 10^{-3}$ \\[-5pt] (9.10)} & \makecell[c]{$6.13\cdot 10^{-3}$ ($1.62\cdot 10^6$, \\[-5pt] $6.83\cdot 10^{3}$)} \\
        \hline
    \end{tabular}
    \caption{System failure probabilities for selected nodes, calculated by the proposed algorithm (c.f. Fig.~\ref{fig:map_ind}) and by the MCS (c.f. Fig.~\ref{fig:map_upd}). Computational time is measured in seconds. MCS is run until c.o.v. reaches 0.01. The numbers in parenthesis are the number of samples (in case of MCS) and computational time.}
    \label{tab:pfs}
\end{table}

\section{Discussion}\label{sec:dis}

\subsection{Connection to previous work}\label{sec:str}

The BRC algorithm builds on investigations made by previous decomposition algorithms. Previous algorithms can be considered special cases of the BRC algorithm in terms of either a system performance function adapted for a branch-and-bound (c.f. Section~\ref{sec:sysFunUseInp}) or a decomposition ordering strategy (c.f. Section~\ref{sec:dec}). For example, in the context of connectivity reliability between two terminals, \cite{LiHe02} proposed a system function that can identify a minimal survival rule from a result of connectivity analysis, which later was improved by \cite{LimSon12} that modified the algorithm to identify the \textit{most likely} survival rules (which is used for the illustrations in Sections~\ref{sec:sysFunUseInp} and \ref{sec:toy}). Similarly, \cite{JanLai08}'s decomposition algorithm for two-terminal maximum flow reliability can be considered a system performance function that can identify a minimal survival rule from a maximum flow analysis result. Meanwhile, building on \cite{LimSon12}, \cite{LeeSon21} improves decomposition efficiency by prioritising component events with a higher centrality measure during decomposition. Though less relevant because of their assumptions on the a-priori availability of minimal rules before decomposition, \cite{LiuEtal21, ForYeh22} similarly proposed effective schemes for decomposition prioritisation based on joint appearances of components in given rules.

The heuristics of the BRC algorithm are developed based on these observations and are found to be effective as supported by the numerical examples in Sections~\ref{sec:toy} and \ref{sec:ex}. However, the BRC algorithm has four distinctive strategies compared to the previous algorithms.

First, to minimise the total number of system performance function runs, it does not run the system performance function directly on the bounds of branches; instead the associated states are inferred from hitherto identified rules. For instance, in Iteration 3 in Fig.~\ref{fig:toy_s3}, $s_{u,2}$, $s_{u,3}$, $s_{u,4}$, and $s_{u,6}$ are all inferred by the rule $\gamma_2 = \left( (x_1^1,x_3^1), 1 \right)$ without additional runs of the system performance function. This is in contrast to previous algorithms that would perform such runs on the bound of each of the branches. We also note that identified minimal rules represent major failure and survival modes, which can be used for planning emergency activities.

Second, since it is impossible to decide the relative importance of components from the beginning (i.e. to decide a decomposition order), the BRC algorithm renews the decomposition process every time it finds a new rule (c.f. Section~\ref{sec:dec}). Thereby, it does not only reduce the number of branches but also identify a component vector state ($x^*$ in Fig.~\ref{fig:brc}) that is likely to lead to a high-probability rule.  This is in contrast to most previous algorithms that do not re-visit branches once they are decomposed. For example, in Iteration 3 in Fig.~\ref{fig:toy_s3}, the BRC algorithm selects $\bm{u}_1$ to check the associated system state, which leads to finding a new failure rule $\gamma_3 = \left( (x_1^0), 0 \right)$.

Third, to remove limitations on the system scale, the BRC algorithm derives not only bounds but also performs hybrid inference (c.f. Section~\ref{sec:stop}). Thereby, even when a system size is too large to be fully decomposed, the BRC algorithm can be employed to effectively reduce the number of system performance function runs, which is usually the most expensive part of an analysis.

Fourth, to be applicable to general system events, the only assumption we introduce for a system performance function is coherency. If advanced understanding of a system event is available, we assure improved efficiency by allowing for an optional output of a system performance function, i.e. an associated rule of a system analysis result (c.f. Section~\ref{sec:sysFunUseInp}). 

Lastly, we note a general property of branch-and-bound algorithms, namely that changes in the probability of component states do not require another round of branch-and-bound. They can be addressed by calculating the probabilities of previously identified branches, which requires only relatively simple arithmetic. By deriving \eqref{eq:upd}, we extend this advantage to the case of hybrid inferences to handle even larger-scale systems. For instance, such probability updating is useful for performing a logic tree analysis during Probabilistic Seismic Hazard Assessment (PSHA) and updating hazard or fragility models \cite{PagaMoneWeat14}.

\subsection{Computational performance}

The BRC algorithm has two primary sources of computational cost: system performance function evaluations and decomposition (c.f. Section~\ref{sec:dec}). The BRC algorithm is designed to minimise the number of system performance function evaluations to handle expensive system functions. As the algorithm finds a new rule with every run of a system function, the required number of system function runs tends to increase linearly with the number of existing rules. This increase is particularly mild when system analysis results can be used to infer an associated survival or failure rule (c.f. Section~\ref{sec:sysFunUseInp}). We note that the number of rules tends to remain limited as real-world systems often have a limited number of survival and failure mechanisms. 

To minimise running a system performance function, the BRC algorithm is thorough in selecting a next component vector state to run a system function on. Hence, it repeatedly performs depth-first decompositions, whose number rapidly increases with the size of a current rule set. While it is difficult to quantify the relationship between the decomposition time and the number of rules precisely, we empirically discuss the results of the example of Section~\ref{sec:ema}. When we compare two system events having 15 rules and 227 rules, the total time for the BRC algorithm takes 2.75 seconds and 16,000 seconds, respectively. These times are mostly spent for decomposition procedures as system performance functions are only run 16 and 228 times, respectively. In other words, in this example, a 15 times increase in the number of rules lead to a 5,807-fold increase in the decomposition time. An increase in decomposition time is more rapid than that in the number of branches. For instance, the two analyses lead to 59 and 25,808 branches, respectively, where the increase is 437-fold.

We note, however, that updates by changed component probabilities do not require repeating decomposition. In the case of the aforementioned examples, updating component probabilities requires 0.0470 and 16.1 seconds, respectively, corresponding to a 342-fold increase. The difference between the two computational times stems from the different numbers of branches. We note that in both cases, the computational times are marginal and are suitable for near real-time analysis.

While efficiency of sampling techniques in general depends on the magnitude of the system failure probability, the efficiency of branch-and-bound algorithms depends on the number of existing rules. These algorithms, including the BRC algorithm, exploit that system events often have a limited number of rules despite their large number of components. This observation is valid especially in systems operated by complex mechanisms such as transportation networks or energy grid systems. In these cases, and when survival rules or failure rules can be inferred from system analysis results, the BRC algorithm can be particularly efficient.

Meanwhile, sparsity of rules may not hold for systems where components play symmetric roles, e.g. series systems, parallel systems, and $k$-out-of-$N$ systems. For instance, in a $k$-out-of-$N$ system, which requires at least $k$ components surviving among $N$ components, there are $\binom{N}{k}$ minimal survival rules and $\binom{N}{N-k+1}$ minimal failure rules, which increase rapidly with $N$. In this case, computational efficiency might be obtained by some other reliability method, e.g. decision diagrams that exploit a symmetry between components \cite{ByuSon21}.

\subsection{Application to multi-state systems}

We note that the assumption of a binary-state system is not restrictive and the BRC algorithm can be applied to multi-state systems. One can perform reliability analysis multiple rounds with different threshold values between survival and failure. In this case, one can minimise the number of running a system performance function by storing analysis results in previous rounds and re-using them for a next round. For instance, if a system functionality is defined by the shortest travel time between an origin and a destination, one can define and analyse multiple thresholds such as 30 minutes, 45 minutes, etc. In order to re-use results from previous thresholds, one can store associated travel times for analysed component vector states. If a previous threshold was 30 minutes, and a system function run identifies that a component vector state $\bm{x}$ leads to a travel time 40 minutes long, this result can be stored for the next analysis with a threshold of 45 minutes so that one can identify $\bm{x}$ as a survival rule without further running an analysis.

\section{Concluding remarks}\label{sec:con}

In this paper, a novel branch-and-bound algorithm is developed for reliability analysis of coherent systems, namely the branch-and-bound for reliability analysis of general coherent systems (BRC) algorithm. The only constraint to apply the BRC algorithm is that a system event should have a coherent definition, i.e. a system state does not improve (worsen) with a worse (better) component vector state, which is valid for many real-world systems. The BRC algorithm improves on previous decomposition algorithms by being applicable to general and/or large-scale systems as well as to multi-state components. To the best of the authors' knowledge, the BRC algorithm is the first branch-and-bound algorithm for reliability analysis that is applicable to \textit{general} coherent systems.

We include a detailed discussion on the distinctive advantages of the BRC algorithm compared to previous algorithms. The numerical investigations demonstrate the significant reduction in the number of system function runs by the BRC algorithm compared to a state-of-the-art algorithm. Further advantages of the BRC algorithm are demonstrated by analysing the Eastern Massachusetts highway (EMA) benchmark network, which consists of 74 nodes and 129 edges. By taking advantage of computational efficiency of the BRC algorithm, we consider each of the nodes a system event and thereby draw maps of failure probabilities over the region.

The BRC algorithm is essentially an algorithm to find minimal rules so that a system event space can be characterised economically. Identified rules by the BRC algorithm have the potential for supporting advanced inference tasks or combining with other reliability methods. For instance, assuming some similarity with unknown rules, identified rules can be used to decide an optimal importance sampling distribution to facilitate hybrid inference. Meanwhile, identified rules can be used to formulate a system performance function as a Boolean function, which allows for employing reliability methods such as hasing-based techniques \cite{SooGocMee20} and probably approximately correct (PAC) technique \cite{ParDueMee19}. Other future research topics include addressing generally coherent but occasionally incoherent systems, performing advanced inference tasks such as component importance measure, and handling with efficiency varying thresholds to obtain the conditional distribution of a system event.

\section*{Data statement}

The Python code for the BRC algorithm is developed as a submodule of the MBNPy toolkit available at \href{https://github.com/jieunbyun/BNS-JT}{https://github.com/jieunbyun/BNS-JT}. The tutorial on using the code is available at \\ \href{https://jieunbyun.github.io/BNS-JT/user-guide/2024/09/15/tutorial_brc}{https://jieunbyun.github.io/BNS-JT/user-guide/2024/09/15/tutorial\_brc}.

\section*{Acknowledgements}

The first author acknowledges support by the Humboldt Foundation through a Humboldt Research Fellowship for Postdoctoral Researchers.

\bibliographystyle{elsarticle-num-names} 
\bibliography{ref_gb}

\appendix
\setcounter{table}{0}
\renewcommand{\thetable}{A\arabic{table}}

\section{Parameters of the EMA benchmark network}

\begin{table}[H]
    \centering
    \begin{tabular}{cc|cc|cc|cc|cc}
    \hline
        Edge & node pair & Edge & node pair & Edge & node pair & Edge & node pair & Edge & node pair \\
        \hline
        $e_{1}$ & $(n_{1}, n_{3})$ & $e_{27}$ & $(n_{14}, n_{22})$ & $e_{53}$ & $(n_{27}, n_{28})$ & $e_{79}$ & $(n_{39}, n_{48})$ & $e_{105}$ & $(n_{51}, n_{52})$ \\
        $e_{2}$ & $(n_{1}, n_{7})$ & $e_{28}$ & $(n_{16}, n_{17})$ & $e_{54}$ & $(n_{27}, n_{33})$ & $e_{80}$ & $(n_{40}, n_{41})$ & $e_{106}$ & $(n_{52}, n_{53})$ \\
        $e_{3}$ & $(n_{1}, n_{9})$ & $e_{29}$ & $(n_{16}, n_{19})$ & $e_{55}$ & $(n_{27}, n_{35})$ & $e_{81}$ & $(n_{40}, n_{48})$ & $e_{107}$ & $(n_{53}, n_{54})$ \\
        $e_{4}$ & $(n_{2}, n_{3})$ & $e_{30}$ & $(n_{16}, n_{22})$ & $e_{56}$ & $(n_{28}, n_{37})$ & $e_{82}$ & $(n_{41}, n_{49})$ & $e_{108}$ & $(n_{55}, n_{57})$ \\
        $e_{5}$ & $(n_{3}, n_{6})$ & $e_{31}$ & $(n_{17}, n_{22})$ & $e_{57}$ & $(n_{29}, n_{41})$ & $e_{83}$ & $(n_{42}, n_{43})$ & $e_{109}$ & $(n_{56}, n_{57})$ \\
        $e_{6}$ & $(n_{3}, n_{7})$ & $e_{32}$ & $(n_{18}, n_{19})$ & $e_{58}$ & $(n_{29}, n_{49})$ & $e_{84}$ & $(n_{42}, n_{45})$ & $e_{110}$ & $(n_{57}, n_{58})$ \\
        $e_{7}$ & $(n_{4}, n_{6})$ & $e_{33}$ & $(n_{18}, n_{21})$ & $e_{59}$ & $(n_{30}, n_{31})$ & $e_{85}$ & $(n_{43}, n_{44})$ & $e_{111}$ & $(n_{57}, n_{59})$ \\
        $e_{8}$ & $(n_{4}, n_{8})$ & $e_{34}$ & $(n_{19}, n_{22})$ & $e_{60}$ & $(n_{30}, n_{60})$ & $e_{86}$ & $(n_{43}, n_{45})$ & $e_{112}$ & $(n_{58}, n_{59})$ \\
        $e_{9}$ & $(n_{5}, n_{10})$ & $e_{35}$ & $(n_{20}, n_{21})$ & $e_{61}$ & $(n_{31}, n_{32})$ & $e_{87}$ & $(n_{44}, n_{46})$ & $e_{113}$ & $(n_{59}, n_{60})$ \\
        $e_{10}$ & $(n_{5}, n_{11})$ & $e_{36}$ & $(n_{20}, n_{30})$ & $e_{62}$ & $(n_{31}, n_{60})$ & $e_{88}$ & $(n_{44}, n_{54})$ & $e_{114}$ & $(n_{59}, n_{72})$ \\
        $e_{11}$ & $(n_{6}, n_{8})$ & $e_{37}$ & $(n_{20}, n_{58})$ & $e_{63}$ & $(n_{32}, n_{33})$ & $e_{89}$ & $(n_{45}, n_{46})$ & $e_{115}$ & $(n_{60}, n_{61})$ \\
        $e_{12}$ & $(n_{6}, n_{13})$ & $e_{38}$ & $(n_{21}, n_{22})$ & $e_{64}$ & $(n_{32}, n_{34})$ & $e_{90}$ & $(n_{45}, n_{47})$ & $e_{116}$ & $(n_{60}, n_{63})$ \\
        $e_{13}$ & $(n_{6}, n_{17})$ & $e_{39}$ & $(n_{21}, n_{23})$ & $e_{65}$ & $(n_{32}, n_{60})$ & $e_{91}$ & $(n_{46}, n_{47})$ & $e_{117}$ & $(n_{60}, n_{65})$ \\
        $e_{14}$ & $(n_{7}, n_{9})$ & $e_{40}$ & $(n_{22}, n_{23})$ & $e_{66}$ & $(n_{33}, n_{34})$ & $e_{92}$ & $(n_{46}, n_{54})$ & $e_{118}$ & $(n_{60}, n_{67})$ \\
        $e_{15}$ & $(n_{7}, n_{13})$ & $e_{41}$ & $(n_{22}, n_{25})$ & $e_{67}$ & $(n_{34}, n_{35})$ & $e_{93}$ & $(n_{47}, n_{48})$ & $e_{119}$ & $(n_{60}, n_{69})$ \\
        $e_{16}$ & $(n_{8}, n_{11})$ & $e_{42}$ & $(n_{22}, n_{28})$ & $e_{68}$ & $(n_{34}, n_{60})$ & $e_{94}$ & $(n_{47}, n_{74})$ & $e_{120}$ & $(n_{60}, n_{71})$ \\
        $e_{17}$ & $(n_{8}, n_{16})$ & $e_{43}$ & $(n_{22}, n_{29})$ & $e_{69}$ & $(n_{35}, n_{36})$ & $e_{95}$ & $(n_{48}, n_{49})$ & $e_{121}$ & $(n_{60}, n_{72})$ \\
        $e_{18}$ & $(n_{9}, n_{12})$ & $e_{44}$ & $(n_{22}, n_{40})$ & $e_{70}$ & $(n_{35}, n_{71})$ & $e_{96}$ & $(n_{48}, n_{51})$ & $e_{122}$ & $(n_{62}, n_{63})$ \\
        $e_{19}$ & $(n_{9}, n_{13})$ & $e_{45}$ & $(n_{23}, n_{24})$ & $e_{71}$ & $(n_{36}, n_{43})$ & $e_{97}$ & $(n_{48}, n_{52})$ & $e_{123}$ & $(n_{63}, n_{65})$ \\
        $e_{20}$ & $(n_{10}, n_{11})$ & $e_{46}$ & $(n_{23}, n_{31})$ & $e_{72}$ & $(n_{36}, n_{44})$ & $e_{98}$ & $(n_{48}, n_{53})$ & $e_{124}$ & $(n_{64}, n_{65})$ \\
        $e_{21}$ & $(n_{10}, n_{18})$ & $e_{47}$ & $(n_{24}, n_{25})$ & $e_{73}$ & $(n_{36}, n_{71})$ & $e_{99}$ & $(n_{48}, n_{54})$ & $e_{125}$ & $(n_{65}, n_{66})$ \\
        $e_{22}$ & $(n_{10}, n_{20})$ & $e_{48}$ & $(n_{24}, n_{26})$ & $e_{74}$ & $(n_{37}, n_{38})$ & $e_{100}$ & $(n_{48}, n_{74})$ & $e_{126}$ & $(n_{67}, n_{68})$ \\
        $e_{23}$ & $(n_{11}, n_{19})$ & $e_{49}$ & $(n_{24}, n_{33})$ & $e_{75}$ & $(n_{37}, n_{42})$ & $e_{101}$ & $(n_{49}, n_{50})$ & $e_{127}$ & $(n_{67}, n_{69})$ \\
        $e_{24}$ & $(n_{13}, n_{14})$ & $e_{50}$ & $(n_{25}, n_{26})$ & $e_{76}$ & $(n_{38}, n_{39})$ & $e_{102}$ & $(n_{49}, n_{73})$ & $e_{128}$ & $(n_{69}, n_{70})$ \\
        $e_{25}$ & $(n_{13}, n_{15})$ & $e_{51}$ & $(n_{26}, n_{27})$ & $e_{77}$ & $(n_{38}, n_{42})$ & $e_{103}$ & $(n_{50}, n_{51})$ & $e_{129}$ & $(n_{69}, n_{71})$ \\
        $e_{26}$ & $(n_{14}, n_{17})$ & $e_{52}$ & $(n_{26}, n_{28})$ & $e_{78}$ & $(n_{39}, n_{40})$ & $e_{104}$ & $(n_{50}, n_{73})$ &   &   \\
        \hline
    \end{tabular}
    \caption{Connectivity of edges (roadways) in EMA highway benchmark network investigated in Section~\ref{sec:ema}}
    \label{tab:ema_edge}
\end{table}

\begin{table}[h!]
    \centering
    \begin{tabular}{c|c c c c c c c c c c c}
    \hline
        \multirow{2}{*}{$M_w$} & \multicolumn{11}{c}{$r_{\text{rm}}$ (km)} \\ \cline{2-12}
        & 1 & 2 & 3 & 5 & 7 & 10 & 20 & 30 & 40 & 50 & 70 \\ \hline
        5.0 & 0.24 & 0.25 & 0.25 & 0.21 & 0.18 & 0.13 & 0.06 & 0.10 & 0.14 & 0.17 & 0.23 \\ \hline
        5.4 & 0.22 & 0.23 & 0.23 & 0.22 & 0.19 & 0.16 & 0.10 & 0.11 & 0.14 & 0.18 & 0.23 \\ \hline
        5.8 & 0.19 & 0.20 & 0.21 & 0.20 & 0.19 & 0.17 & 0.14 & 0.15 & 0.18 & 0.21 & 0.27 \\ \hline
        6.2 & 0.15 & 0.15 & 0.16 & 0.16 & 0.17 & 0.16 & 0.14 & 0.15 & 0.18 & 0.21 & 0.27 \\ \hline
        6.6 & 0.12 & 0.13 & 0.13 & 0.13 & 0.14 & 0.14 & 0.13 & 0.15 & 0.18 & 0.21 & 0.27 \\ \hline
        7.0 & 0.13 & 0.13 & 0.13 & 0.13 & 0.14 & 0.14 & 0.14 & 0.15 & 0.17 & 0.20 & 0.27 \\ \hline
        7.4 & 0.14 & 0.13 & 0.13 & 0.14 & 0.15 & 0.16 & 0.16 & 0.16 & 0.18 & 0.21 & 0.27 \\ \hline
        7.8 & 0.15 & 0.15 & 0.14 & 0.15 & 0.15 & 0.17 & 0.18 & 0.18 & 0.20 & 0.22 & 0.28 \\ \hline
        8.2 & 0.17 & 0.17 & 0.16 & 0.16 & 0.17 & 0.18 & 0.20 & 0.21 & 0.22 & 0.23 & 0.29 \\ \hline
    \end{tabular}
    \caption{Epistemic standard deviation $\sigma_{\ln Y}^e$ with $Y$ as Sa at period 1.0 sec. The values are excerpted from Table 1 in Appendix of \cite{Cam03}.}
    \label{tab:EpiStd}
\end{table}

\begin{table}[h!]
    \centering
    \begin{tabular}{c|c|c|c}
        \hline
        \makecell{HAZUS \\[-5pt] classification} & Description & \makecell{Mean ($R$) and dispersion ($\beta$) of Sa (g)\\[-5pt]at $T=1.0$ sec (extensive damage state)} & Road IDs \\
        \hline
        HWB2 & \makecell{Non CA; seismic;\\[-5pt]major bridge} & 1.10 and 0.6 & \makecell{$e_{27}$, $e_{30}$, $e_{31}$, $e_{34}$, $e_{38}$, $e_{40}$,\\[-5pt]$e_{41}$, $e_{42}$, $e_{43}$, $e_{44}$, $e_{58}$}\\
        \hline
        HWB4 & \makecell{Non CA, seismic;\\[-5pt]single span} & 1.20 and 0.6 & \makecell{$e_{11}$, $e_{20}$, $e_{24}$, $e_{28}$, $e_{32}$, $e_{45}$,\\[-5pt]$e_{47}$, $e_{48}$, $e_{50}$, $e_{51}$, $e_{52}$, $e_{53}$,\\[-5pt]$e_{59}$, $e_{61}$, $e_{64}$, $e_{70}$, $e_{74}$, $e_{76}$,\\[-5pt]$e_{78}$, $e_{80}$, $e_{85}$, $e_{86}$, $e_{87}$, $e_{89}$,\\[-5pt]$e_{90}$, $e_{91}$, $e_{123}$, $e_{126}$} \\
        \hline
        HWB11 & \makecell{Non CA; seismic;\\[-5pt]continuous concrete} & 1.10 and 0.6 & All other roads \\
        \hline
    \end{tabular}
    \caption{Fragility curve classifications of roadways ($e_1.\ldots,e_{129}$) assumed for analysis and associated fragility parameters.}
    \label{tab:ema_fra}
\end{table}

\end{document}